\documentclass[preprint,12pt]{elsarticle}
\usepackage{blindtext}
\usepackage[a4paper, total={6in, 8in}]{geometry}
\usepackage{array,multirow,graphicx}
\usepackage{float}
\usepackage{amsfonts, amsmath, amssymb, amsthm}
\usepackage{comment}
\usepackage{xcolor}
\usepackage{subcaption}
\usepackage{siunitx}
\usepackage{arydshln}
\usepackage{bbm}
\usepackage{lineno}
\usepackage{setspace}

\newtheorem{problem}{Problem}
\newtheorem*{problemOL}{Problem}

\newtheorem*{ocproblem}{Optimal Control Problem}

\newtheorem*{kolm}{KOL-m}

\newtheorem*{kolp}{KOL-$\partial$}

\newtheorem*{remark}{Remark}

\DeclareMathOperator*{\argmin}{arg\,min}
\newcommand{\NNML}{NN-ModL }
\newcommand{\NNMLL}{NN-ModL}

\title{Learning epidemic trajectories through Kernel Operator Learning: from modelling to optimal control}
\affiliation[mox]{organization={MOX, Dipartimento di Matematica, Politecnico di Milano},
      addressline={Piazza Leonardo da Vinci 32}, 
      city={Milano},
      postcode={20133}, 
      country={Italy}}
\author[mox]{Giovanni Ziarelli}
\author[mox]{Nicola Parolini}
\author[mox]{Marco Verani}
\date{April 2024}

\makeatletter
\def\ps@pprintTitle{%
  \let\@oddhead\@empty
  \let\@evenhead\@empty
  \let\@oddfoot\@empty
  \let\@evenfoot\@oddfoot
}
\makeatother

\begin{document}
\begin{abstract}
Since infectious pathogens start spreading into a susceptible population, mathematical models can provide policy makers with reliable forecasts and scenario analyses, which can be concretely implemented or solely consulted.
In these complex epidemiological scenarios, 
machine learning architectures can play an important role, since they directly reconstruct data-driven models circumventing the specific modelling choices and the parameter calibration, typical of classical compartmental models.
In this work, we discuss the efficacy of Kernel Operator Learning (KOL)
to reconstruct population dynamics during epidemic outbreaks, where the transmission rate is ruled by an input strategy.
In particular, we introduce two surrogate models, named KOL-m and KOL-$\partial$, which reconstruct in two different ways the evolution of the epidemics.
Moreover, we evaluate the generalization performances of the two approaches with different kernels, including the Neural Tangent Kernels, and compare them with a classical neural network model learning method.
Employing synthetic but semi-realistic data, we show how the two introduced approaches are suitable for realizing fast and robust forecasts and scenario analyses, and how these approaches are competitive for determining optimal intervention strategies with respect to specific performance measures.
\textbf{Highlights:}
\begin{itemize}
    \item We formalize the novel Kernel Operator Learning (KOL) framework in the context of epidemic problems;
    \item We build two KOL strategies for estimating susceptible, infectious and recovered individuals given an input level of reduction of transmission rate; 
    \item We numerically test the two approaches on synthetic data for assessing their reliability in making scenario-analyses;
    \item We compare the solutions of two optimal control problems on standard compartmental models with respect to the solutions retrieved with KOL approaches.
\end{itemize}

\end{abstract}
\maketitle

\section{Introduction}
The recent global SARS-CoV-2 pandemic has underlined the paramount importance of developing mathematical models and numerical schemes for predictions and forecasts of epidemic illnesses: from the perspective of policy-makers, it is often useful to dispose of qualitative and quantitative results for making scenario analyses and forecasts; from the social point of view sharing information about possible outcomes can be beneficial in order to increase social awareness and public knowledge on the illness current spreading and its future evolutions.
A typical approach relies on traditional compartmental mathematical frameworks, where specific modelling choices and parameters embody the different features that characterize the epidemic spreading, the virological effects of the illness and the impact of different pharmaceutical interventions.
However, in presence of new epidemic outbreaks many key-features could still be unknown or difficult to isolate in clinical trials, and consequently the illness itself could be difficult to be completely described by the classical compartmental models.
Indeed, clinical symptoms of different illnesses are multifaceted and strictly depend on the origin of the pathogenic microbial agent responsible of the disease, which can be bacterial, parasitic, fungal, viral or originated by prions, \textit{i.e.} other kinds of toxic proteins, and on the pathway through which the illness naturally diffuses \cite{martcheva2015introduction, brauer2019mathematical}.
Moreover, in order to accurately describe the disease through compartmental models, it is fundamental to account for possible DNA or RNA mutations from the wildtype strain in long-term outbreaks,
as well as for possible preventive measures and control, including vaccination, treatments, prophylaxis, quarantine, isolation or other measures minimizing social activities (like the use of face-masks, compulsory home-schooling, different levels of lockdowns).
In these highly complex and rapidly changing scenarios, the efficacy of compartmental models for making fast scenario-analyses may be severely limited by the delicate and sometimes ad-hoc parameter calibration process that becomes even harder if one aims at embodying age-dependency or other geographical features, see, \textit{e.g.}, \cite{giordano2020modelling, parolini2022modelling, bertuzzo2020geography}. 

Alongside with scenario analysis, recent upcoming epidemic events have shown the importance of disposing of computational tools measuring the impact of 
pharmaceutical resources \cite{gozzi2022anatomy} and other Non-Pharmaceutical Interventions (NPIs) \cite{he2023combining} 
so to guide policy-makers in choosing how to intervene limiting the social and economic burden.
From the mathematical perspective, we can leverage on the versatility of optimal control theory in order to derive useful quantitative and qualitative guidelines for minimizing the amount of infectious or deceased individuals \cite{ziarelli2023optimized, lemaitre2022optimal}, the total incidence of the spreading disease \cite{britton2023optimal}, or the amount of contacts and, consequently, of cases \cite{dimarco2022optimal}.
Other problems which have been further investigated, analytically and numerically, are more delicate from the mathematical viewpoint such as the minimization of epidemic peaks \cite{morris2021optimal} or the minimization of the eradication time \cite{bolzoni2017time}.

In view of the above discussion, it is of paramount importance to provide the society with mathematical tools able to output computationally cheap and reliable scenario analyses, so to compare different prevention measures and solve optimal control problems.
Among recently developed mathematical frameworks, a prominent position is covered by Operator Learning, which, roughly speaking, deal with the development and application of algorithms designed to build approximation of operators starting from a given set of input/output pairs. 
In the family of Operator Learning tools, an increasing attention has been devoted to the so called Deep Operator Networks (DeepONets) introduced in \cite{lu2021learning}.
Since then, other ML algorithms exploiting deep neural networks have been developed, such as PINNDeepONets \cite{goswami2022physics} and Fourier Neural Operators \cite{li2020fourier}.
More recently, in \cite{batlle2024kernel} Kernel Operator Learning (KOL) has been proposed as a competitive alternative to the previous approaches in terms of cost-accuracy trade-off and the capability of matching (sometimes outperforming) in several benchmark problems the generalization properties of learning methods based on neural networks.
Moreover, the simple closed formula of the learnt kernel operator makes KOL very attractive.

In this work we propose two different numerical approaches based on Kernel Operator Learning (KOL) \cite{batlle2024kernel}, that starting from epidemic data provide surrogate models describing the dependency of different stages of the illness on a given control function representing NPIs. To the best of our knowledge, these approaches are new in the epidemic context. Moreover, the numerical experiments contained in the sequel of the paper show that the presented approaches can be efficiently employed for making fast scenario analyses and solving optimization tasks, since they can be rapidly trained directly from data, circumventing delicate calibration phases. 
For the sake of simplicity and to present the main features of our approaches, we work with synthetic data generated by standard epidemic differential models governed by control functions modelling NPIs (or any other effect reducing the transmissibility of the illness), which we assume to be given.

As pointed out in \cite{batlle2024kernel}, the performance of KOL depends on the choice of the kernel employed to perform the regression tasks, see also \cite{genton2001classes}.
Very recently, Neural Tangent Kernels (NTKs) have been introduced in a series of pioneering works (cf. \cite{lee2020finite, jacot2018neural, arora2019exact}) that opened the door to a very intensive research activity.
Loosely speaking, NTK arise from the connection with infinite-width neural networks and, within the infinite-width limit, the generalization properties of the neural network could be described by the generalization properties of the associated NTK.
In view of this connection and encouraged by the results presented in \cite{lee2020finite, adlam2023kernel}, in this paper we employ NTKs in the construction of our kernel operator approaches and we validate its efficacy through a wide campaign of numerical tests.

The paper is organized as follows.
In Section \ref{sec:SECMAT}, we synthetically illustrate and derive the mathematical formulation of our KOL approaches and we briefly recall the compartmental models that we adopt for generating the numerical data to test the methods.
In Section \ref{sec:SECRES} we present and discuss different numerical tests for assessing the generalization properties of our approached together with their efficacy in providing solutions to optimization tasks.

\section{Materials and Methods}
\label{sec:SECMAT}
This section is organized in three parts: Section \ref{sec:kolEpi} is devoted to introducing Kernel Operator Learning in the epidemic context, whilst Section \ref{sub:overEpi} briefly recalls the standard compartmental epidemic models that will be employed to produce the synthetic data feeding the KOL.
In Section \ref{app:preliminaries} we gather some preliminary considerations about the numerical technicalities for retrieving both KOL regressors.
\subsection{KOL and epidemic modelling: Basic principles}
\label{sec:kolEpi}
\begin{figure}[t]
  \centering
  \includegraphics[width=\textwidth]{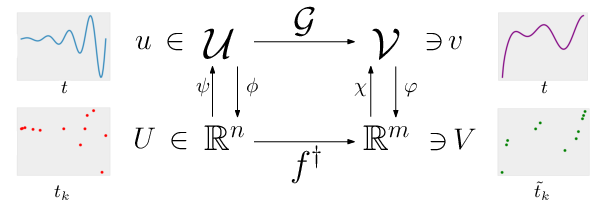}
  \caption{Kernel Operator Learning (KOL) diagram. Operating among RKHS justifies the reduction of the problem to learning the behaviour of the vector-value function $f^\dagger$, operating between input/output observations (cf. Appendix \ref{app:kolDer}).}
  \label{fig:diagram}
\end{figure}
We briefly summarize some of the principal results on Kernel Operator Learning (KOL) contained in  \cite{batlle2024kernel}, that will be instrumental to build surrogate epidemic models and to solve optimal control problems.

To this aim, let us consider two possibly infinite-dimensional Hilbert spaces $\mathcal{U},\mathcal{V}$ and assume there exists an unknown operator mapping between the two spaces, \textit{i.e.}
\begin{equation}
  \mathcal{G}: \mathcal{U} \rightarrow \mathcal{V}.
\end{equation}
Roughly speaking, the goal of operator learning is to approximate $\mathcal{G}$ based on pairs of input/output that are accessible through finite dimensional linear measurements, as formalized by the following:
\begin{problemOL}
  Let $\{u_i, v_i\}_{i=1}^N$ be $N$ samples in $\mathcal{U} \times \mathcal{V}$, \textit{i.e.}
  \begin{equation}
    \mathcal{G}(u_i) = v_i, \; \mathrm{with} \, i = 1,2 \hdots N.
  \end{equation}
  Moreover, define the (bounded and linear) observation operators $\phi : \mathcal{U} \rightarrow \mathbb{R}^n$ and $\varphi : \mathcal{V} \rightarrow \mathbb{R}^m$ acting on the input and the target functions, respectively.
  The goal of operator learning is the approximation of the operator $\mathcal{G}$ based on the observation input/output pairs $\{ \phi(u_i), \varphi(v_i) \}_{i=1}^N$.
\end{problemOL}
In the rest of the section, to ease the reading, we restrict the presentation to the scalar case, while the extension to the vector-valued case, relevant for the numerical examples in the following sections, can be straightforwardly obtained from the scalar one by resorting to the theory of vector-valued Reproducing Kernel Hilbert Spaces (RKHS) in \cite{carmeli2006vector}.
Therefore, let us consider $\mathcal{U}$ and $\mathcal{V}$ as functional spaces made of scalar functions, where the dependant variable is named $t$ and it is assumed to vary in the interval $D \subset \mathbb{R}$.
In this scenario, a standard choice for the observation operators $\phi$ and $\varphi$ consists in considering the pointwise evaluation at specific collocation points $\{t_k\}_{k=1}^m$ and $\{\Tilde{t}_k\}_{k=1}^n$ respectively, which in general can be different (see Figure \ref{fig:diagram}).
However, for simplicity, in this work we consider the same collocation points, \textit{i.e.}
\begin{equation}
  \begin{split}
  \phi : u \rightarrow U := (u({t}_1), u({t}_2) \hdots u({t}_n))^T \in \mathbb{R}^n,\\
  \varphi : v \rightarrow V := (v(t_1), v(t_2) \hdots v(t_n))^T \in \mathbb{R}^n.
  \end{split}
\end{equation}
Then, assuming we are given the training dataset $\{ U_j, V_j\}_{j=1}^N$ where, consistently with our notation, $U_j := (u_j({t}_1), u_j({t}_2) \hdots u_j({t}_n))^T \in \mathbb{R}^n$ and 
$V_j := (v_j(t_1), v_j(t_2) \hdots v_j(t_n))^T \in \mathbb{R}^n$, we aim at constructing an approximation $\bar{\mathcal{G}}$ of $\mathcal{G}$.
More specifically, following \cite{batlle2024kernel}, endowing $\mathcal{U}$ and $\mathcal{V}$ with a RKHS structure and using kernel regression to identify the maps $\psi$ and $\chi$ (cf. Figure \ref{fig:diagram} and see Appendix \ref{app:kolDer} for more details), the approximated operator can be written explicitly in closed form as
\begin{equation}
  \bar{\mathcal{G}}(u)(t) = K(t, T) K(T,T)^{-1} \left ( \sum_{j=1}^N S(\phi(u), U_j) \alpha_j \right ),
  \label{eq:KOLGeneral}
\end{equation}
where $K$ is the kernel function induced by the RKHS structure of $\mathcal{V}$, the vector $T = [t_1, t_2 \hdots t_n]^T \in \mathbb{R}^n$ contains the collocation points, and $S:\mathbb{R}^n \times \mathbb{R}^n \rightarrow \mathbb{R}$ is a properly chosen vector-valued kernel.
Moreover, $K(\cdot, T) : D \rightarrow \mathbb{R}^n$ is a row vector such that $K(t,T)_i = K(t, t_i)$, and $K(T,T)$ is an $n \times n$ matrix such that $K(T,T)_{ij} = K(t_i, t_j)$.
Parameters $\{\alpha_j\}_{j=1}^N \in \mathbb{R}^n$ are the kernel regression parameters over the input/output training pairs.
Since we consider pointwise observation operators, $K(t,T)$ is the composite linear interpolant of nodes in $T$ evaluated at the desired time $t$.
We refer to \cite{batlle2024kernel} for a complete discussion of the convergence properties of this approach and for the formal derivation of the learnt operator in a more general functional framework.

Although the above framework is quite general and can be applied to a variety of problems, here, given the goal of the paper, we embody it in the context of reconstructing processes ruled by dynamical systems steered by control variables.
More precisely, we consider the following general scalar dynamical system ruled by the control function $u$
\begin{equation}
\begin{cases}
    \dot{{v}}(t) = {F}({v}(t), {u}(t)), \; \forall t \in (0, t^*]\\
    {v}(0) = {v}_0
\end{cases}
\label{eq:dynSys}
\end{equation}
where $F: \mathbb{R} \times \mathbb{R} \rightarrow \mathbb{R}$ is assumed to be sufficiently smooth to guarantee the well-posedness of the problem \cite{fleming2012deterministic}.
Clearly, we can associate to \eqref{eq:dynSys} the mapping $\mathcal{G}$ that given the control $u$ returns the solution $v$.
This mathematical framework fits a typical epidemic context, where it is often necessary to take into account the impact of different NPIs acting in reducing the contact rate among vulnerable people, such as, for instance, partial or total lockdowns, the mandatory use of face-masks in public spaces or the isolation of people showing mild symptoms which could be linked to the illness of concern.

Now, we are ready to introduce our strategies.
The first strategy, named KOL-m, consists in directly determining an approximation of the solution map $\cal{G}$ from $u(t)$ to $v(t)$, as described in the following.
\begin{kolm}
  Let $\mathcal{U} = \{ u \in L^2(0,t^*) \}$ and $\mathcal{V} = C^0([0, t^*])$. 
  Given the input-output data trajectories $\{ (\hat{u}_k, \hat{v}_k)\}_{k=1}^{N}$, 
  where each $\hat{u}_k \in \mathcal{U}$ is a control function and each $\hat{v}_k \in \mathcal{V}$ is the associated (known) state vector function, the Kernel Operator Learning Map method (KOL-m) builds 
  $\mathcal{\bar{G}}_m: \mathcal{U} \rightarrow \mathcal{V}$ according to \eqref{eq:KOLGeneral}.
\end{kolm}
Hence, the quantity $ \mathcal{\bar{G}}_m\left ( u^*\right)(t)$ represents an approximation to the solution $v(t)$ of \eqref{eq:dynSys} with control variable $u^*(t)$.

Let us now describe the second approach, namely KOL-$\partial$, which approximates the operator that given the prescribed control returns the derivative of the solution to \eqref{eq:dynSys}.
\begin{kolp}
  Let $\mathcal{U} = \{ u \in L^2(0,t^*) \}$ and $\hat{\mathcal{V}} = L^1(0,t^*)$. 
  Given the input-output data trajectories $\{ (\hat{u}_k, \dot{\hat{v}}_k)\}_{k =1}^{N}$, where each $\hat{u}_k \in \mathcal{U}$ is a control function and each $\dot{\hat{v}}_k \in \hat{\mathcal{V}}$ is the associated (known) derivative of the state vector function, the Kernel Operator Learning Derivative method (KOL-$\partial$) builds $ \mathcal{\bar{G}_\partial}: \mathcal{U} \rightarrow \hat{\mathcal{V}} $ 
  according to \eqref{eq:KOLGeneral}.
  \end{kolp}
  
  Therefore, the quantity $ \mathcal{\bar{G}}_\partial\left ( u^*\right)(t)$ represents an approximation to $\dot{v}(t)$, \textit{i.e.} the derivative of the solution to \eqref{eq:dynSys} with control variable $u^*(t)$.
  Specifically, we have
  \begin{equation}
    v(t) = v_0 + \int_{0}^t \mathcal{\bar{G}}_\partial ( u) (\tau) d\tau \; \in \mathcal{V}
    \; \mathrm{where} \; \mathcal{\bar{G}}_\partial: \mathcal{U} \rightarrow \hat{\mathcal{V}} \; \mathrm{satisfies} \; 
    \eqref{eq:KOLGeneral},
  \end{equation}
  or alternatively
  \begin{equation}
  \begin{cases}
  \begin{split}
  \dot{v}(t) & = K(t, T) K(T,T)^{-1} \left ( \sum_{j=1}^N S(\phi(u), U_j) \alpha_j \right )\\
  &= K(t, T) K(T,T)^{-1} \left ( \sum_{j=1}^N S(\phi(u), U_j) \begin{pmatrix}
    [\textbf{S}(\textbf{U}, \textbf{U})^{-1} \hat{\textbf{V}}_{\cdot, 1}]_j\\
    [\textbf{S}(\textbf{U}, \textbf{U})^{-1} \hat{\textbf{V}}_{\cdot, 2}]_j\\
    \vdots\\
    [\textbf{S}(\textbf{U}, \textbf{U})^{-1} \hat{\textbf{V}}_{\cdot, n}]_j\\  
  \end{pmatrix} \right ) \;\;\; \forall t \in (0,t^*]
  \end{split}\\
  v(0) = v_0,
  \end{cases}
  \label{eq:KOLparComplete}
  \end{equation}
  where $\hat{\textbf{V}}_{\cdot, k} = [[\varphi(\dot{\hat{v}}_1)]_k, [\varphi(\dot{\hat{v}}_2)]_k \hdots [\varphi(\dot{\hat{v}}_N)]_k]^T, \; \forall k=1,2\hdots n$ is the vector of the evaluations at the $k$-th point of each output, $\textbf{U}$ is a vector collecting all the $\{U_j\}_j^N$ and $\textbf{S}(\textbf{U}, \textbf{U}) \in \mathbb{R}^{N \times N}$ can be defined as $[\textbf{S}(\textbf{U}, \textbf{U})]_{ij} = S(U_i, U_j)$ (cf. Appendix \ref{app:kolDer}).
  In this case we assume that an accurate approximation of the derivative of each state $( \{ v_i \}_{i=1}^N)$ is observable, hence the sequence $ \{ \varphi(v_i) \}_{i=1}^N$ is available (cf. Figure \ref{fig:diagram}).

We conclude this section with a methodological remark on the application of the KOL-m approach to the reconstruction of compartments when dealing with epidemic problems or other problems where the reconstructed operator has to preserve positivity.
This modified version of KOL-m will be employed to derive the numerical results presented in the subsequent sections.
\begin{remark}[On the positivity preserving property of KOL-m]
The solutions to the epidemic differential problems are intrinsically positive, since they represent positive fractions of the population corresponding to different states with respect to the illness spreading. However, there is no reason why KOL-m should preserve the positivity of the prediction, even in presence of positive data $\{\hat{\textbf{x}}_k\}$. 
For this reason, in order to enforce the positivity of the prediction of the learnt operator $\mathcal{\bar{G}}_m$, we proceed as follows: given $\{ (\hat{u}_k, \hat{\textbf{x}}_k)\}_{k \leq N}$, with positive $\hat{\textbf{x}}_k$, we apply KOL-m to the modified input-output dataset $\{ (\hat{u}_k, \sqrt{\hat{\textbf{x}}_k})\}_{k \leq N}$ so to obtain the intermediate operator $\mathcal{\widetilde{G}}_m$. Then, we set $\mathcal{\bar{G}}_m=\mathcal{\widetilde{G}}^2_m$.
The efficiency of this approach is showed in Section \ref{sec:KOL-numerics}.
\end{remark}
\subsection{Overview of compartmental models}
\label{sub:overEpi}
\begin{figure}[h!]
  \centering
  \includegraphics[width=\textwidth]{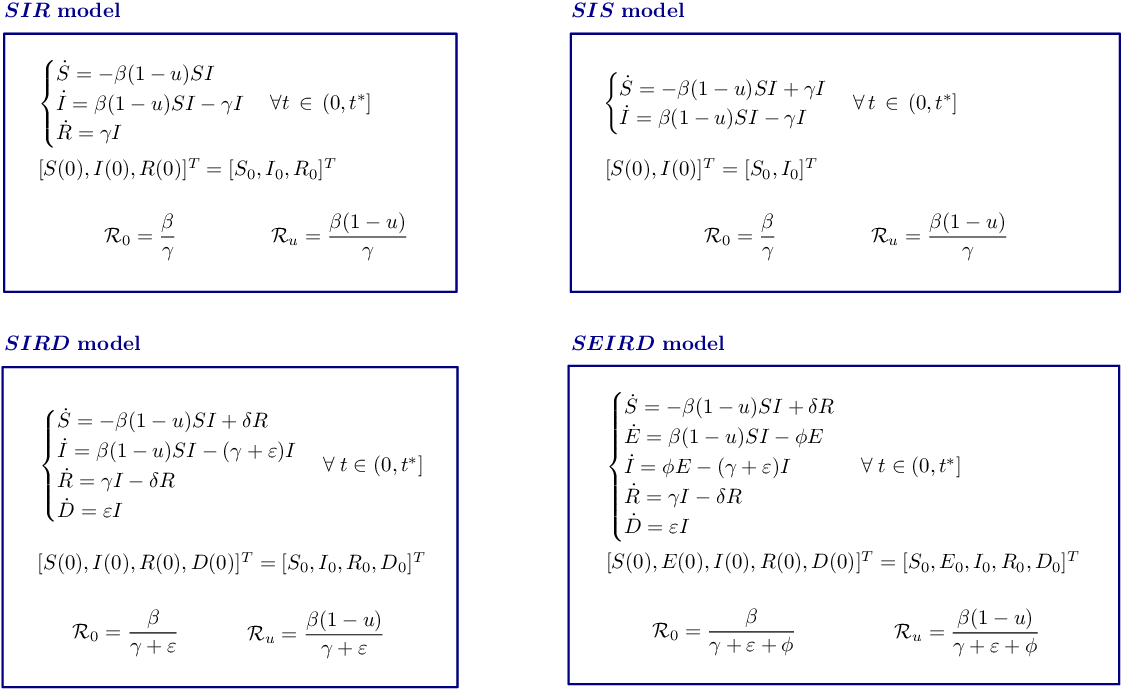}
  \caption{$SIR$, $SIS$, $SIRD$ and $SEIRD$ compartmental models.
  Each model is provided with its respective basic reproduction number ($\mathcal{R}_0$) and the reproduction number depending on the control ($\mathcal{R}_u$).}
  \label{fig:epimods}
\end{figure}

In order to test KOL-m and KOL-$\partial$ we generate synthetic data employing classical compartmental models and use the latter to validate the accuracy of the surrogate models.
More precisely, for simplicity and without loss of generality, we restrict our attention to the four classical compartmental models reported in Figure \ref{fig:epimods}, where the control variable $u$ represents the instantaneous reduction of the basic transmission rate $\beta$ dictated by the virological and transmissibility properties of the illness by reducing the contact rate.
We assume that $u_{lb} \leq u(t) \leq u_{ub}, \; \forall \, t \in [0,t^*] $, where the upper bound value of the control function $u_{ub} = 1$ models total lockdown, instead the lower bound $u_{lb} = 0$ stands for null NPIs.
For suitable choices of the function $\textbf{F}: \mathbb{R}^d \times \mathbb{R} \rightarrow \mathbb{R}^d$, all those models (and other high dimensional disease-specific compartmental models, \textit{e.g.} \cite{parolini2021suihter, marziano2021effect, bertuzzo2020geography}) can be written in the following general form:
\begin{equation}
\begin{cases}
    \dot{\boldsymbol{x}}(t) = \textbf{F}(\boldsymbol{x}(t), u(t)), \; \forall t \in (0, t^*]\\
    \boldsymbol{x}(0) = \boldsymbol{x}_0.
\end{cases}
\label{eq:dynSys2}
\end{equation}
Clearly, the solution of \eqref{eq:dynSys2} can be written in terms of the operator $\mathcal{G}$, that given the input control $u(t):[0, t^*]\rightarrow\mathbb{R}$ returns the state $\boldsymbol{x}(t): [0,t^*] \rightarrow \mathbb{R}^d$.
For each model, $\forall \, t \, \in \, [0, t^*]$ the evolution function $\textbf{F}$ is sufficiently smooth in order to guarantee existence and uniqueness of the solution of Cauchy problem \eqref{eq:dynSys2} when $u$ admits at most a countable amount of jump discontinuities (see, \textit{e.g.} \cite{yong2012stochastic}).
For a more comprehensive review on epidemiological models, we refer, \textit{e.g.}, to \cite{martcheva2015introduction, brauer2019mathematical}.

\subsection{Numerical aspects of KOL-m and KOL-$\partial$}
\label{app:preliminaries}

In this section we collect some useful information on the computational aspects involved in the training process of KOLs.
For the purpose of this work we assume that the space of our discrete counterpart of the control functions and the corresponding discrete vectors coming from the evaluation of the solution are of the same dimension, \textit{i.e.} $m=n$ in Figure \ref{fig:diagram}.
It is worth to remark that, in the epidemic context, it is natural to consider, as observation operators for the state and the control, the pointwise evaluation, which is coherent with the way open access datasets for epidemic events are often organized and real-time forecasts are delivered (see, \textit{e.g.} \cite{sherratt2023predictive}).

For what concerns the generation of synthetic data, we employ $N$ different control functions to be observed via the observation operators, together with the $N$ associated compartmental trajectories (vector in $\mathbb{R}^d$) obtained by numerically solving one of the differential systems in Figure \ref{fig:epimods} using Explicit Euler with time step $dt$. 

From the computational perspective, once the scalar kernel $S(\cdot, \cdot)$ between the two discrete vector spaces has been chosen, we need to solve $n$ linear systems of dimension $N \times N$ where the matrix is $\textbf{S}(\textbf{U},\textbf{U})$ (cf. Appendix \ref{app:ntk}).
For this purpose, we compute the Cholesky Factorization of the matrix and solve the systems with the standard substitution methods (cf. \cite{meanti2020kernel} for more advanced strategies).
Finally, in solving the regression problem we add a regularization term, with penalty parameter equal to $1e-10$.



\section{Results and Discussion}
\label{sec:SECRES}
\label{sec:KOL-numerics}
In this section we present a wide campaign of numerical tests with the aim of: (a) understanding the impact of the choice of different kernels in KOL-m and KOL-$\partial$ (cf. Section \ref{sub:ntk}); (b) comparing our KOL methods with a popular neural-network based model learning method (cf. Section \ref{sub:compNN}); (c) assessing the robustness of the introduced approaches for solving two optimization tasks (cf. Section \ref{sec:kolContr}).

\subsection{On the choice of the Kernel}
\label{sub:ntk}
One crucial ingredient, driving the approximation and generalization properties of the corresponding KOL method, is the choice of the scalar kernel $S$ in \eqref{eq:GammaK}.
Which kernel function is optimal for kernel regression is still an open debate and depends on the specific application.
For instance, there exists some innovative approaches which learn the kernels by simulating data driven dynamical systems enabling scalability of Kernel Regression \cite{owhadi2019kernel}.
Here, we consider the following popular choices for $S$. 

\begin{figure}[H]
  \centering
  \begin{subfigure}[b]{0.35\textwidth}
    \centering
    \includegraphics[width=\textwidth]{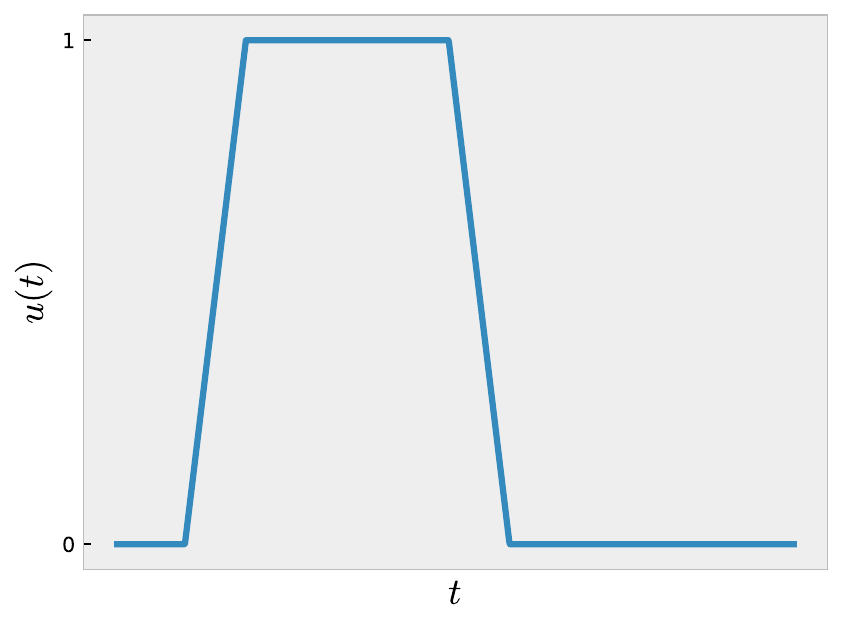}
    \caption{}
    \label{fig:u1}
  \end{subfigure}
  \begin{subfigure}[b]{0.35\textwidth}
    \centering
    \includegraphics[width=\textwidth]{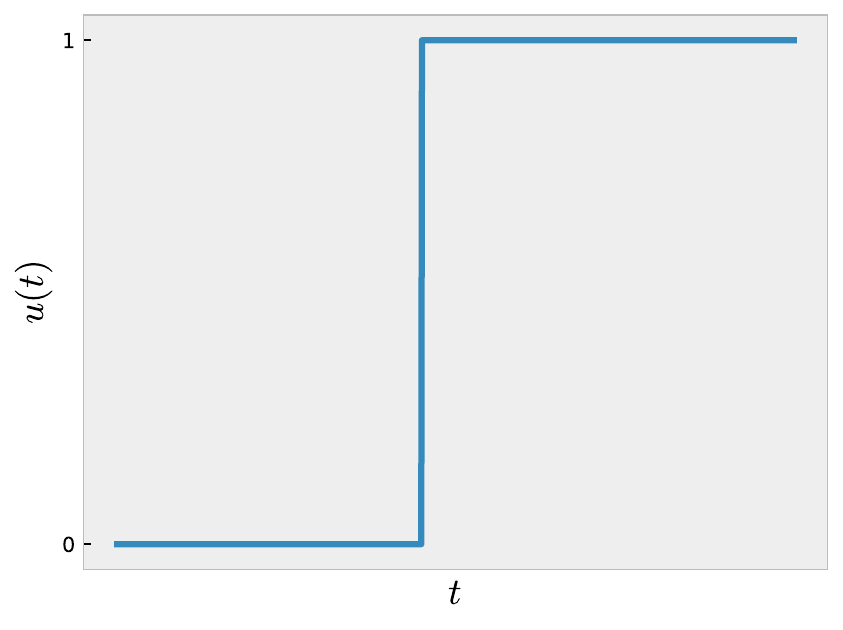}
    \caption{}
    \label{fig:u2}
  \end{subfigure}
  \\
  \begin{subfigure}[b]{0.35\textwidth}
    \centering
    \includegraphics[width=\textwidth]{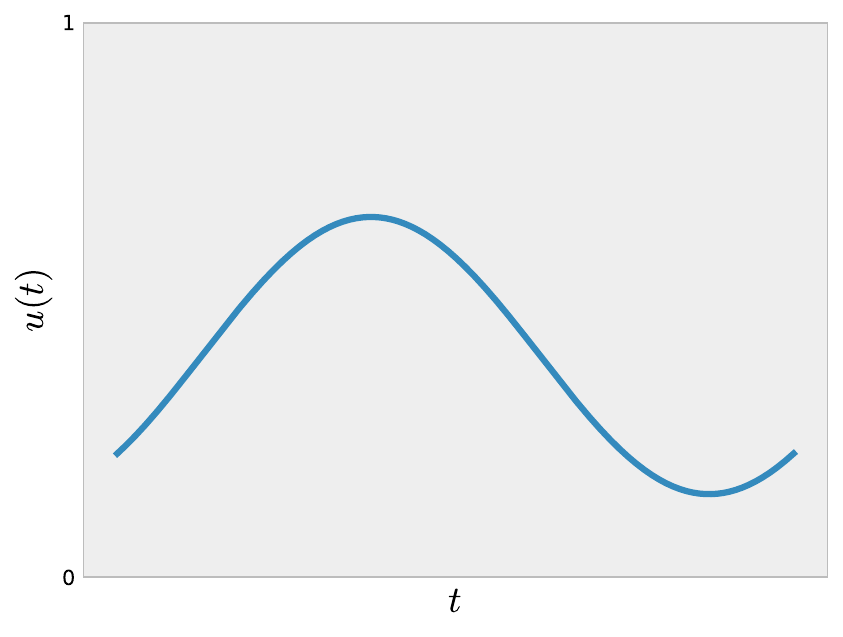}
    \caption{}
    \label{fig:u3}
  \end{subfigure}
  \begin{subfigure}[b]{0.35\textwidth}
    \centering
    \includegraphics[width=\textwidth]{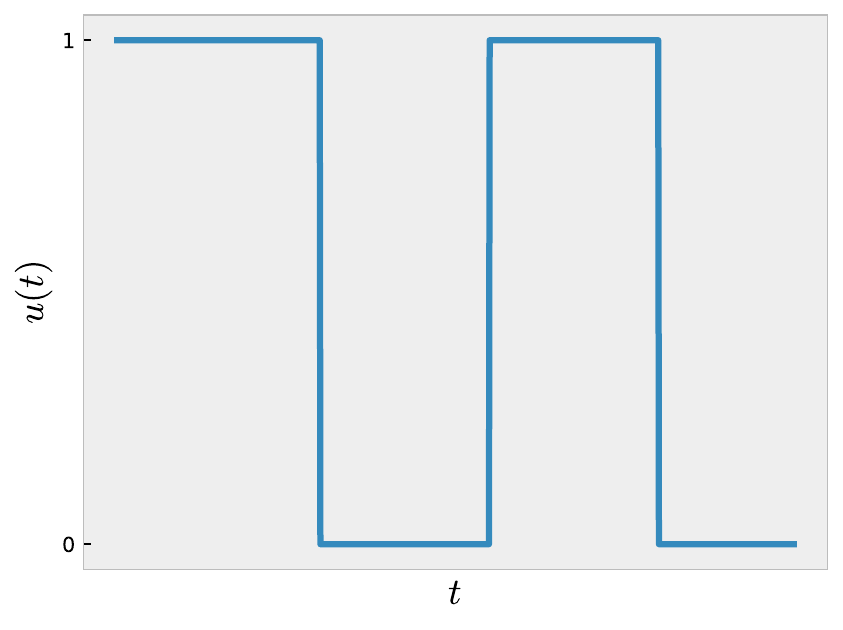}
    \caption{}
    \label{fig:u4}
  \end{subfigure}
  \caption{Examples of control functions employed to generate the training-testing dataset for the KOL methods.}
  \label{fig:controls}
\end{figure}

\begin{itemize}
  \item \textbf{Linear Kernel}: $S(U_1, U_2) = U_1^T U_2$.
  This kernel evaluates the alignment in the $n$-dimensional space between input vectors;
  \item \textbf{Matérn Kernel}: $S(U_1,U_2) = \dfrac{2^{1-\nu}}{\Gamma(\nu)}\left ( \sqrt{2 \nu} \dfrac{\| U_1 - U_2\|_2}{\rho}\right)^{\nu} K_{\nu} \left ( \sqrt{2 \nu} \dfrac{\| U_1 - U_2\|_2}{\rho}\right) $, where $\nu > 0$ controls the smoothness of the kernel function, $ \rho $ is a characteristic length scale, $\Gamma$ is the Gamma function and $K_{\nu}$ is the modified second kind Bessel function of order $\nu$. 
  This family of kernels is stationary and isotropic, given the Euclidean distance between input points.
  If $\nu \rightarrow \infty$ one obtains the Gaussian kernel.
  It is often used for image analysis \cite{de2020shape} and other machine learning regression tasks \cite{gulian2019machine};
  \item \textbf{RBF kernel}: $S(U_1, U_2) = e^{-\frac{\| U_1 - U_2\|_2^2}{2 \sigma^2}}$, known also as Gaussian kernel.
  It has the interpretation of a similarity measure since it is bounded in $[0,1]$ and decreases as long as the distance between points increases;
  \item \textbf{Rational Quadratic kernel}: $S(U_1, U_2) = \left ( 1 + \dfrac{\| U_1 - U_2\|_2^2}{2 \alpha l^2}\right )^{-\alpha}$, with $l,\alpha > 0$. 
  It can be regarded as the infinite sum of different RBF kernels with different length scales;
  \item \textbf{Neural Tangent Kernel (NTK)}: Given a neural network regressor $f(x; \theta)$ of depth $d_{nn}$, width $l_{nn}$ and activation function $\sigma_{nn}$, with $\theta$ denoting the vector collecting all weights and biases, we define the family of finite-width Neural Tangent Kernels $ \{ S \}_{\tau > 0}: \mathbb{R}^n \times \mathbb{R}^n \rightarrow \mathbb{R} $ (cf. \ref{app:ntk}) as
  \begin{equation}
    S_\tau(x_i,x_j) := \langle \partial_{\theta} f(x_j; \theta(\tau)) ,\, \partial_{\theta} f(x_i; \theta(\tau)) \rangle,
  \end{equation}
  where $\tau$ represents a fictitious iteration time. 
  It has been proven that, if the initialization of the weights follows the so-called NTK initialization \cite{jacot2018neural}, in the infinite-width limit each element in $\{ S_\tau\}_{\tau}$ converges in probability to a stationary kernel independently on $\tau$, \textit{i.e.}
  \begin{equation}
    S_{\tau}(x_i, x_j) \underset{\mathbb{P}}{\rightarrow} S(x_i, x_j), \; \forall \, \tau > 0, \; \forall \, x_i,x_j.
  \end{equation}
  Hence, the family of NTKs strictly depends on two parameters: activation function and depth of the associated neural network.
  
  In this paper, with NTK we refer to the infinite-width limit kernel function.
  We rely on the library presented in \cite{novak2019neural} for an efficient and user-friendly computation of the considered NTKs.
  In practice, given the network hyperparameters, we compute the scalar NTK $S: \mathbb{R}^d \times \mathbb{R}^d \rightarrow \mathbb{R}$ and then evaluate it at desired couples of points.
  More details about the derivation of NTKs associated to neural networks can be found in Appendix \ref{app:ntk}.

\end{itemize}
For each of the above kernels, we evaluate the generalization properties of KOL-m and KOL-$\partial$ trained on the epidemic $SIS$, $SIR$ and $SEIRD$ trajectories generated by control functions chosen among the following families (see Figure \ref{fig:controls}) which are representative of possible interventions that policy makers can implement.

For training and testing KOL-m and KOL-$\partial$, we define a mixed dataset constituted by functions belonging to four distinct functional families (see Figure \ref{fig:controls}):
\begin{enumerate}
  \item \textbf{Linear Pulse} (Figure \ref{fig:u1}): \begin{equation}
    u(t) = \begin{cases}
      u_0, & t \leq t_0,\\[10pt]
      \dfrac{3 (u_1 - u_0)}{\Delta t} (t - t_0) + u_0 & t_0 < t \leq t_0 + \frac{\Delta t}{3},\\
      u_1 & t_0 + \frac{\Delta t}{3} < t \leq t_0 + \frac{4 \Delta t}{3}, \\[10pt]
      \dfrac{3 (u_0 - u_1)}{\Delta t} (t - t_0 - \frac{4 \Delta t}{3}) + u_1 & t_0 + \frac{4 \Delta t}{3} < t \leq t_0 + \frac{5 \Delta t}{3},\\
      u_0 & t > t_0 + \frac{5 \Delta t}{3},
    \end{cases}
  \end{equation}
  depending on 4 degrees of freedom (dof): $u_0, \, u_1 \in [0,1]$, $t_0 \in [0,t^*]$ and $\Delta t \in [0, \frac{t^*}{3}]$;
  \item \textbf{Step function} (Figure \ref{fig:u2}): 
  \begin{equation}
    u(t) = \begin{cases}
      u_0 & t \leq t_0,\\
      u_1 & t > t_0,
    \end{cases}
  \end{equation}
  with 3 dofs: $u_0, \, u_1 \in [0,1]$, $t_0 \in [0,t^*]$;
  \item \textbf{Continuous Seasonality} (Figure \ref{fig:u3}):
  \begin{equation}
    u(t) = \dfrac{u_0}{2} \left ( 1 + \dfrac{1}{2} \cos{ \left ( \dfrac{2 \pi t}{t^*} + \dfrac{\Delta t}{t^*} \dfrac{\pi}{2} \right )} \right ), 
  \end{equation}
  with 2 dofs: $u_0 \in [0,1]$ and $\Delta t \in [0, \frac{t^*}{3}]$;
  \item \textbf{Double step} (Figure \ref{fig:u4}):
  \begin{equation}
    u(t) = \begin{cases}
      u_0 & t \leq t_0,\\
      u_1 & t_0 < t \leq t_0 + \frac{\Delta t}{2},\\
      u_0 & t_0 + \frac{\Delta t}{2} < t \leq t_0 + \Delta t,\\
      u_1 & t > t_0 + \Delta t,\\
    \end{cases}
  \end{equation}
  with 3 dofs: $u_0, \, u_1 \in [0,1]$, $t_0 \in [0,t^*]$.
\end{enumerate}
To test the generalization properties we proceed as follows.
For each given input control function $\{u^{n}\}_{n}$, $n=1, 2, \hdots, N_p$ with $N_p$ number of test points, we generate the output data $\{\boldsymbol{x}^{i,n}\}_{i,n}$, $i=1,2,\hdots,N_c$ with $N_c$ the number of compartments, by numerically solving the specific systems of ODEs ($SIS$, $SIR$, $SIRD$ and $SEIRD$) up to a final time $t^* = 100$, with a discretization step $dt = 1$.
The pairs control/output data form the training dataset.
We define the prediction relative error as
\begin{equation}
  p_{\mathrm{err}} = \displaystyle\dfrac{1}{N} \sum_{n=1}^{N}\displaystyle\sum_{i=1}^{N_c} \dfrac{\| \boldsymbol{x}_p^{i,n} - \boldsymbol{x}^{i,n}\|_2}{\| \boldsymbol{x}^{i,n}\|_2},
  \label{eq:perror}
\end{equation}
where $\{\boldsymbol{x}_p^{i,n}\}$ are the predictions generated by the KOL methods.
Clearly, the prediction samples belong to a different batch with respect to the training input data.

In all cases, the value of the basic reproduction number is set at $\mathcal{R}_0 = 4$, and the recovery rate $\gamma = 0.05$.
The epidemic event starts with 1\% of infected (or, respectively, exposed) individuals in all cases.
The rest of the population is assumed to be susceptible to the disease at the initial time.
For the $SEIRD$ model, the remaining parameters have been chosen as $\delta = 0.4, \, \varepsilon = \phi = 0.05$.
The transmission rate $\beta$ without control is, therefore, computed starting from the definition of $\mathcal{R}_0$ of each model in Figure \ref{fig:epimods}.

For the comparison of the different generalization properties of the KOL acting on different kernels, the training dataset is built employing 500 control functions equally distributed among the above 4 classes of control functions, where the dofs are randomly sampled following a uniform distribution in the respective interval of definition of each parameter.
As for the prediction samples, we test the KOL approaches over a dataset built employing 100 samples equally distributed among the same families of control functions, but with different values for the dofs.

\begin{figure}[H]
  \centering
  \includegraphics[width=0.8\textwidth]{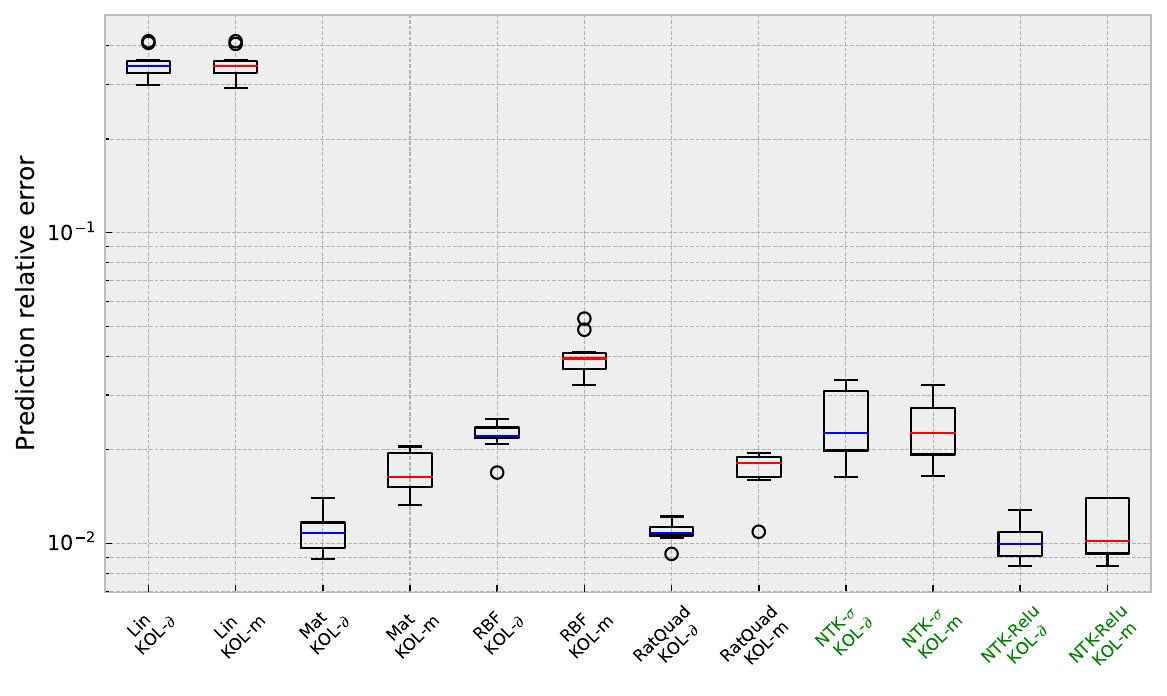}
  \caption{($SIS$ model) Comparison of the prediction errors over 100 control functions with KOL methods trained on 20 batches of size 500, where the control functions are chosen in the mixed training dataset. 
  Bullet points represent outliers whose prediction error is outside the 1.5 IQR (length of the whiskers) of the set of simulations.}
  \label{fig:boxplots_kernels_sis}
\end{figure}

\begin{figure}[H]
  \centering
  \includegraphics[width=0.8\textwidth]{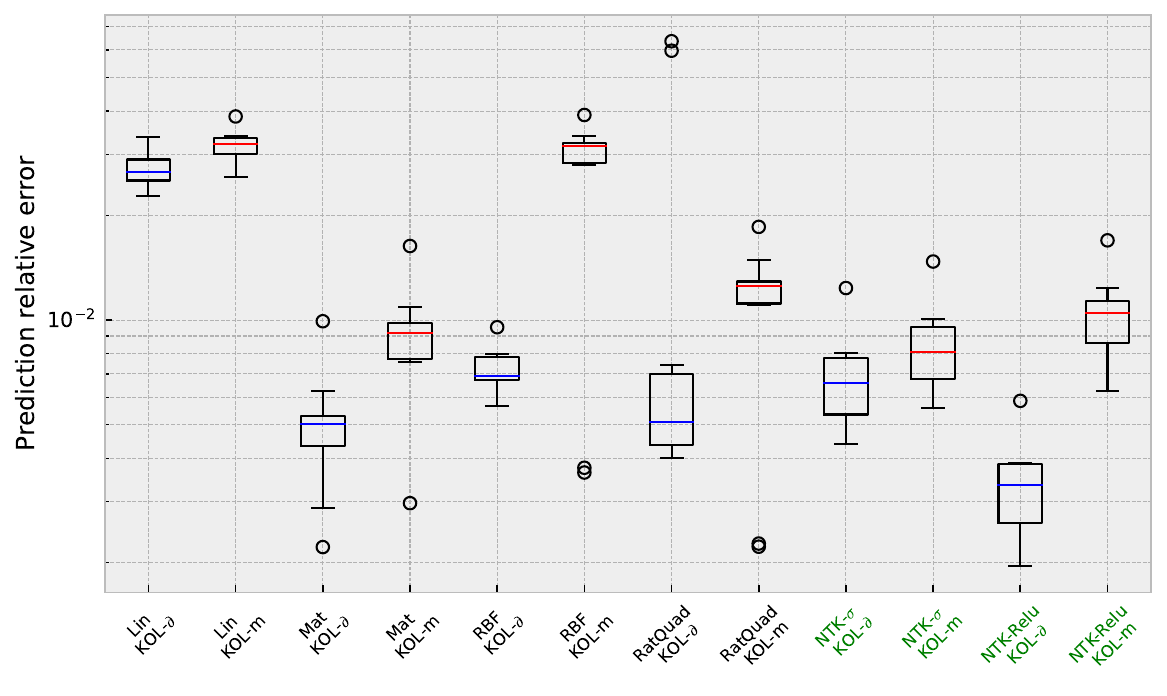}
  \caption{($SIR$ model) Comparison of the prediction errors over 100 control functions with KOL methods trained on 20 batches of size 500, where the control functions are chosen in the mixed training dataset.
  Bullet points represent outliers whose prediction error is outside the 1.5 IQR (length of the whiskers) of the set of simulations.}
  \label{fig:boxplots_kernels_sir}
\end{figure}

\begin{figure}[H]
  \centering
  \includegraphics[width=0.8\textwidth]{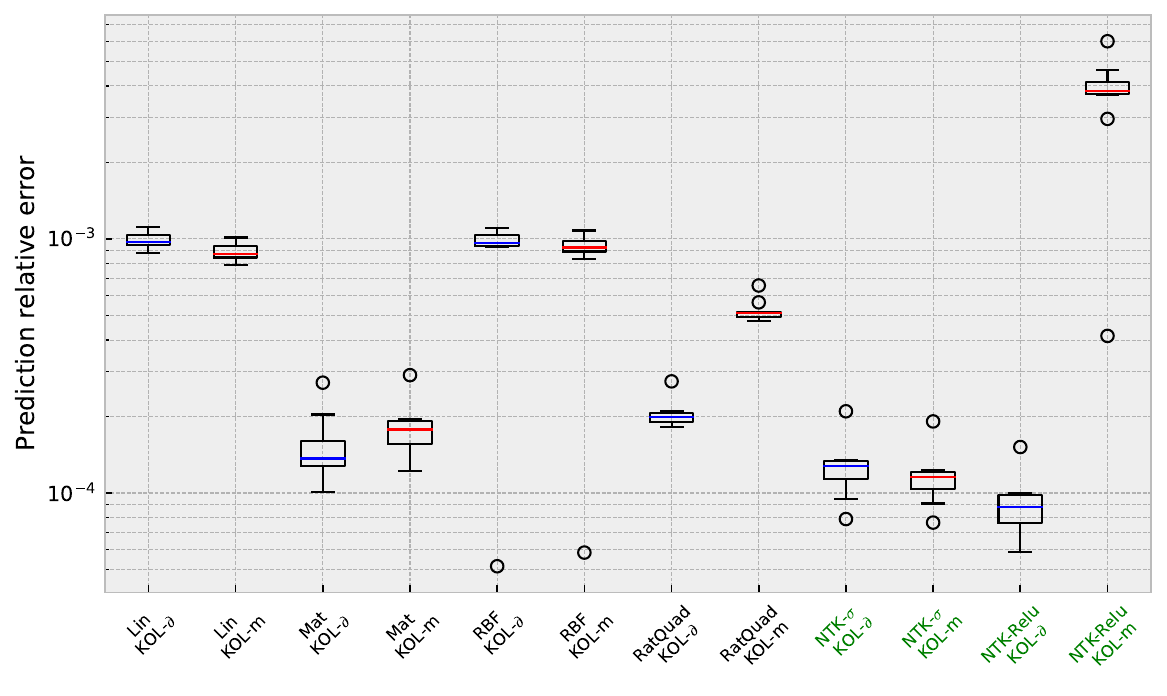}
  \caption{($SEIRD$ model) Comparison of the prediction errors over 100 control functions with KOL methods trained on 20 batches of size 500, where the control functions are chosen in the mixed training dataset.
  Bullet points represent outliers whose prediction error is outside the 1.5 IQR (length of the whiskers) of the set of simulations.}
  \label{fig:boxplots_kernels_seird}
\end{figure}

Results in Figures \ref{fig:boxplots_kernels_sis}-\ref{fig:boxplots_kernels_seird} show the prediction relative error of the different approaches.
In particular, for each kernel we trained each KOL-m and KOL-$\partial$ with 20 different training datasets of size 500 as detailed above and we represented the results in the boxplot form.
We perform the same error analysis for different compartmental models: $SIS$ (Figure \ref{fig:boxplots_kernels_sis}), $SIR$ (Figure \ref{fig:boxplots_kernels_sir}) and $SEIRD$ (Figure \ref{fig:boxplots_kernels_seird}).
The boxplots corresponding to the blue median lines are the one referring to the KOL-$\partial$ approach, whilst the red ones correspond to the KOL-m approach.
For those kernels depending on parameters to be tuned, such as RBF, Matérn and Rational Quadratic, we underdid a sensitivity analysis showing similar results exploring the space of parameters (Matérn: $\nu \in [0,0.1], \, \rho \in [0.1,1]$, RBF: $\sigma \in [0, 0.1]$, Rational Quadratic: $\alpha \in [0.01, 0.1], \, l \in [0.01, 1]$).
For what concerns the NTKs, which depend on the choice of the depth $d_{nn}$ and on the activation function $\sigma_{nn}$ of the associated neural network, we perform the same sensitivity analysis considering both Neural Tangent Kernels with Rectified Linear Unit ($ReLu$) and Sigmoidal activation functions for different depths ($d_{nn} \in \{2,3,4 \}$).
In the boxplots in Figures \ref{fig:boxplots_kernels_sis}-\ref{fig:boxplots_kernels_seird} we reported the results of the NTKs corresponding to two-layer networks, which showed to be a suitable trade-off between computational complexity and generalization properties.

We notice that for the $SIS$ and the $SIR$ model, the prediction errors are higher when compared to the one obtained for the $SEIRD$ model.
This seems to indicate that KOL methods generalize with lower errors for growing complexity.
For the $SIS$ model, the operator approximated via the Linear Kernel does not succeed in surrogating the epidemic dynamics.
On the other hand, median errors of KOL methods with $ReLu$-based NTK achieve the lower results, granting errors always lower than 2\% across the test samples.
Instead, for the $SEIRD$ and $SIR$ model, all the KOLs with the different kernels seem to be good proxies in terms of generalization properties.
Specifically, the KOL-$\partial$ methods have better approximation properties, with a particular mention to the NTK-$Relu$ based.
With the 5 dimensional model $SEIRD$, among the KOL-m approaches the lowest median error with less IQR is attained by the NTK-sigmoidal one.

In virtue of these results, from now on, we select the NTK-sigmoidal kernel for the KOL-m method, and NTK-$ReLu$ kernel for the KOL-$\partial$ method. 

\subsection{Comparison between KOL and a Model Learning based on neural networks}
\label{sub:compNN}
In this section, we enlight some important features of our KOL approaches that make them reliable and competitive tools for model learning. To this aim, we compare KOL-m and KOL-$\partial$ with a popular and paradigmatic machine learning approach based on neural networks, namely the approach described in \cite{regazzoni2019machine}, in the sequel called \NNMLL, and shortly described in what follows.

Let $\mathcal{F} = \{ \boldsymbol{f}_{nn}: \mathbb{R}^{d} \times \mathbb{R} \rightarrow \mathbb{R}^{d} \}$ be the space of feed-forward neural network functions.
The \NNML problem reads as follows:
\begin{problem}[\NNML] Solve the constrained optimization problem
  \begin{equation}
      \min_{\boldsymbol{f}_{nn} \, \in \mathcal{F}} \frac{1}{2}\sum_{j=1}^{N} \int_0^{t^*} |\textbf{x}_j(t) - \textbf{x}_{nn,j}(t)|^2 dt,
  \end{equation}
  s.t.
  \begin{equation}
    \begin{cases}
      \dot{\textbf{x}}_{nn,j}(t) = \textbf{f}_{nn}(\textbf{x}_{nn,j}(t), u_j(t)) & t \in (0,t^*], \, j = 1,2\hdots N,\\
      \textbf{x}_{nn,j}(0) = \textbf{x}_0,
    \end{cases}
    \label{eq:nnml}
  \end{equation}
  \label{prob:nnml}
  where $\boldsymbol{f}_{nn} \, \in \, \mathcal{F}$ and $N$ is the size of the training set.
\end{problem}
Hence, \NNML 
reconstructs the map from the control $u(t)$ to the state $\boldsymbol{x}_{nn}(t)$.
The problem is recast as a discrete optimization problem where the optimizing variables are weights and biases associated to $\boldsymbol{f}_{nn}$.
\NNML has shown valuable performances for building reduced models of cardiac electromechanics \cite{regazzoni2022machine}, or for deriving a data-driven model for active force generation in cardiomyocites \cite{regazzoni2020machine}.
As illustrated in \cite{regazzoni2019machine}, the discretized finite-dimensional version of Problem \ref{prob:nnml} is equivalent to a nonlinear least-square problem, which can been solved employing the Levenberg-Marquardt iterative method \cite{nocedal1999numerical}.

In the sequel, we compare our KOL approaches with \NNML in the epidemic context, in terms of: (a) the wall-clock time employed for training the operator learning; (b) the generalization error in the testing phase.
The comparison has been drawn by accounting for progressively larger training set sizes for both KOL approaches and \NNML in order to approximate the population dynamics generated synthetically by $SIR$, $SIRD$ and $SEIRD$ models.
In particular, we trained \NNML with different training sizes and different number of maximum iterations.
We performed a sensitivity analysis on the choice of the number of neurons per layers in the \NNML case, and concluded, following the \textit{Occam's razor principle of parsimony}, to go for a shallow neural network of 6 neurons.

All the results have been obtained by executing the code in parallel on a 8-core machine Intel i7.
In Table \ref{tab:wallclock} we summarize the outputs of the comparison between the \NNML and the two KOL approaches, where in each case we considered $dt = 0.05, \; t^* = 5, \; \mathcal{R}_0 = 2$, $\gamma = 0.05$, $\delta = 0.4$ and $ \varepsilon = \phi = 0.05$. 
The initial conditions correspond to those of Section \ref{sub:ntk}.
The comparison has been done by considering 50, 100 or 200 iterations of the optimization scheme for the \NNMLL, trained with 25, 50 and 100 elements uniformly sampled in the mixed training dataset.
Both KOL approaches have also been trained on larger datasets constituted by 200, 500, 800, 1000, 2000, 5000, and 10000 elements. 
The prediction error has been evaluated according to \eqref{eq:perror} over the same test set of 100 samples constituted by four batches of dimension 25, each belonging to the four functional families in the mixed dataset.
The wall-clock time of the running processes has been measured in seconds.

For each epidemic model, the wall-clock time taken for the training stage of \NNML is of orders of magnitude higher than the one needed for KOL approaches, even when considering solely 50 iterations.
The \NNML approach ends with a prediction error close to $10^{-2}$ in the three cases after a total wall-clock time ranging from $10^2$ (best scenario) and $10^3$ seconds (worst scenario).
Instead, for all models in less than 8 seconds the KOL-m approach trained with 500 samples reaches a generalization error which is always lower than the corresponding ones reached by \NNMLL.
In addition, for all three models, up to 2000-sample training sets, both KOL approaches still employ less wall-clock time than \NNMLL. 
In particular, in the case of $SEIRD$ model, the prediction error is 3 orders of magnitude lower than the one of \NNMLL.
The gain in terms of wall-clock time is more evident for $SEIRD$ model with size 5000 against \NNML with size 100.

As expected, further increasing the amount of training samples makes the prediction errors continuously decrease though with a slower rate with respect to the training size.
For what concerns the KOL-$\partial$ schemes, they achieve lower prediction errors with respect to the other approaches even with the baseline size of 100 input functions, except for $SIRD$.
Moreover, we deduce that in both KOL-m and KOL-$\partial$ approaches, the prediction errors are not affected by the increasing in dimensions of the differential problem going from the 3 dimensional $SIR$ to the 5 dimensional $SEIRD$.
Even the computational times surprisingly do not undergo significant changes across the different models.

In conclusion, both KOL approaches return more accurate predictions in a lower computational time with respect to the \NNMLL, trained for different number of iterations.
Between the two KOL approaches, KOL-$\partial$ is preferable in terms of accuracy across different epidemic differential models. 
In general, the KOL approaches turn out to possess a desirable property, \textit{i.e.} the ability of returning reliable predictions in a fast way, which is extremely useful in the epidemic context, for instance when scenario analyses or optimization tasks are required. 
For the sake of fairness it should be mentioned that, regarding \NNML, it is theoretically possible to find a neural-network function able to build an approximation of the target state with any desired accuracy.
However, the proof of the existence of such architecture is not constructive.

\begin{table}[h]
  \centering
  \resizebox{\textwidth}{!}{
  \begin{tabular}{ccc|cc|cc|cc}
    \hline
    \multicolumn{3}{c|}{} & \multicolumn{2}{c|}{$SIR$} & \multicolumn{2}{c|}{$SIRD$} & \multicolumn{2}{c}{$SEIRD$} \\
    \hline
    & \textbf{Iterations} & \textbf{Size} & \textbf{Error} & \textbf{Time [s]} & \textbf{Error} & \textbf{Time [s]} & \textbf{Error} & \textbf{Time [s]}\\
    \hline
    \hline
    \NNML & 50 & 25 & $\SI{1.5e-2}{}$ & 279 & $\SI{1.0e-2}{}$ & 380 & $\SI{1.1e-2}{}$ & 383\\
    & 100 & 25 & $\SI{1.5e-2}{}$ & 478 & $\SI{1.0e-2}{}$ & 615 & $\SI{1.1e-2}{}$ & 755 \\
    & 200 & 25 & $\SI{1.5e-2}{}$ & 1186 & $\SI{1.0e-2}{}$ & 1139 & $\SI{1.1e-2}{}$ & 1673\\
    & 50 & 50 & $\SI{1.2e-2}{}$ & 576 & $\SI{1.1e-2}{}$ & 559 &  $\SI{1.1e-2}{}$ & 875\\
    & 100 & 50 & $\SI{1.2e-2}{}$ & 1363 & $\SI{1.1e-2}{}$ & 1456 & $\SI{1.1e-2}{}$ & 1476\\
    & 200 & 50 & $\SI{1.2e-2}{}$ & 3448 & $\SI{1.1e-2}{}$ & 2116 & $\SI{1.1e-2}{}$ & 3046 \\
    & 50 & 100 & $\SI{8.8e-3}{}$ & 893 & $\SI{9.2e-3}{}$ & 1394 & $\SI{9.5e-3}{}$ & 1364 \\
    & 100 & 100 & $\SI{8.8e-3}{}$ & 1672 & $\SI{9.2e-3}{}$ & 2801 & $\SI{9.5e-3}{}$ & 3262\\
    & 200 & 100 & $\SI{8.8e-3}{}$ & 1725 & $\SI{9.2e-3}{}$ & 4511 & $\SI{9.4e-3}{}$ & 6668 \\
    \hline
    KOL-m & & 25 & $\SI{8.3e-2}{}$ & 0.425 & $\SI{8.2e-2}{}$ & 0.385 & $\SI{7.5e-2}{}$ & 0.209 \\
    & & 50 & $\SI{5.4e-2}{}$ & 0.521 & $\SI{5.5e-2}{}$ & 0.486 & $\SI{4.8e-2}{}$ & 0.479 \\
    & & 100 & $\SI{1.6e-2}{}$ & 0.725 & $\SI{1.5e-2}{}$ & 0.731& $\SI{4.9e-3}{}$ & 0.755\\
    & & 200 & $\SI{1.2e-2}{}$ & 1.47& $\SI{1.1e-2}{}$ & 1.492 & $\SI{2.6e-3}{}$ & 1.46\\
    & & 500 & $\SI{5.6e-3}{}$ & 7.60& $\SI{5.4e-3}{}$ & 7.09 & $\SI{9.1e-5}{}$ & 6.85\\
    & & 800 & $\SI{5.5e-3}{}$ & 17.8& $\SI{5.3e-3}{}$ & 18.8 & $\SI{7.0e-5}{}$& 20.7\\
    & & 1000 & $\SI{4.7e-3}{}$ & 32.5 & $\SI{4.5e-3}{}$ & 25.7 & $\SI{6.4e-5}{}$ & 29.4\\
    & & 2000 & $\SI{3.0e-3}{}$ & 116 & $\SI{2.9e-3}{}$ & 94.5 & $\SI{3.4e-5}{}$ & 95.8\\
    & & 5000 & $\SI{2.3e-3}{}$ & 572& $\SI{2.2e-3}{}$ & 568 & $\SI{1.3e-5}{}$ & 559\\
    & & 10000 & $\SI{1.5e-3}{}$ & 2269& $\SI{1.4e-3}{}$ & 2265 & $\SI{8.4e-6}{}$ & 2288 \\
    
     \hline
    KOL-$\partial$ & & 25 & $\SI{7.0e-3}{}$ & 0.420 & $\SI{1.2e-2}{}$ & 1.27 & $\SI{2.6e-3}{}$ & 0.501\\
    & & 50 & $\SI{4.2e-3}{}$ & 0.524 & $\SI{8.6e-3}{}$ & 0.872 & $\SI{1.2e-3}{}$ & 0.856\\
    & & 100 & $\SI{2.5e-3}{}$ & 0.744 & $\SI{4.9e-3}{}$ & 1.17 & $\SI{3.6e-4}{}$ & 1.11 \\
    & & 200 & $\SI{1.5e-3}{}$ & 1.54& $\SI{3.5e-3}{}$ & 2.14 & $\SI{2.3e-4}{}$ & 2.05\\
    & & 500 & $\SI{8.0e-4}{}$ & 6.32& $\SI{1.9e-3}{}$ & 7.85 & $\SI{9.7e-5}{}$ & 8.26 \\
    & & 800 & $\SI{5.9e-4}{}$ & 19.3& $\SI{2.0e-3}{}$ & 18.9 & $\SI{6.2e-5}{}$& 18.1\\
    & & 1000 & $\SI{5.8e-4}{}$ & 27.3 & $\SI{1.8e-3}{}$ & 25.3 & $\SI{5.2e-5}{}$& 26.4\\
    & & 2000 & $\SI{1.2e-4}{}$ & 93.8& $\SI{1.5e-3}{}$ & 90.5 & $\SI{3.9e-5}{}$ & 94.7\\
    & & 5000 & $\SI{1.2e-4}{}$ & 877 & $\SI{1.1e-3}{}$ & 546 & $\SI{1.9e-5}{}$ & 888 \\
    & & 10000 & $\SI{1.2e-4}{}$ & 2206 & $\SI{8.5e-4}{}$ & 2156 & $\SI{1.3e-5}{}$ & 3315\\
     \hline
     \hline
  \end{tabular}}
  \caption{Wall-clock time comparison for the $SIR$, $SIRD$ and $SEIRD$ models.}
\label{tab:wallclock}
\end{table}

\subsection{KOL and epidemic control}
\label{sec:kolContr}
In this section we exploit KOL-m and KOL-$\partial$ for the solution of paradigmatic optimal control problems, which are meaningful for making fast and reliable scenario-analyses in the context of epidemic control. More precisely, in Section \ref{sub:ocgroppi} we are interested in the minimization of the eradication time prescribed a given epidemic threshold, while in Section \ref{sub:ocQuad} we tackle the minimization of the total amount of infected.

\subsubsection{Optimal control for estimating the minimum eradication time}
\label{sub:ocgroppi}

In the sequel, we focus on the minimization of the eradication time prescribed a given epidemic, building upon the theoretical results given in \cite{bolzoni2017time}.
The presence of rigorous mathematical results on the existence of an optimal solution, makes the eradication problem a desirable and mathematically solid benchmark to study the reliability of our KOL-based approach.
We formulate the problem by considering the standard deterministic $SIR$ model.

\begin{problem}[Minimum eradication time]
Let $\mathcal{U} = \{ 0 \leq u(t) \leq u_{max} \, \forall \, t \in [0, t^*] \}$ be the space of admissible controls. Solve
\begin{equation*}
  \min_{u \in \mathcal{U}} \mathcal{J}_{t_e}(u) = \int_{0}^{t_e(u)} 1 \,dt,
\end{equation*}
subject to the state problem
  \begin{equation*}
    \begin{cases}
    \dot{S} = -\beta (1 - u) S I, \\
    \dot{I} = \beta (1-u) S I - \gamma I,\\
    \dot{R} = \gamma I, 
    \end{cases}
  \end{equation*}
$\forall t \, \in (0,t^*]$ given $[S(0), I(0), R(0)]^T = [S_0, I_0, R_0]^T$, where $t_e(u)$ is the eradication time associated to the control $u$.
\label{prob:groppiOC}
\end{problem}

Given a threshold $\eta > 0$ of infected individuals, the eradication time is defined as the first time when infectious individuals are below the given threshold value, \textit{i.e.}
\begin{equation}
  t_e \in (0,t^*] \; \mathrm{s.t.}\; I(t) > \eta \, \forall \, 0 < t < t_e, \, I(t_e) = \eta. 
\end{equation}
It has been proven in \cite{bolzoni2017time} that the space of admissible optimal controls can be restricted to
\begin{equation}
  \mathcal{A} = \{ u:[0,t^*] \rightarrow \{0,u_{max} \}, \; u(t) = u_{max}\, H(t - \tau) \},
\end{equation}
where $\tau$ is the starting intervention time depending on the maximum level of intervention $u_{max}$, and $H(\cdot)$ is the Heaviside step function.
Defining the \textit{maximum control reproduction number} as 
\begin{equation}
  \mathcal{R}_{u_{max}} = \frac{\beta (1 - u_{max})}{\gamma},
\end{equation}
it is possible to find an optimal solution with non-trivial switching time, \textit{i.e.} $\tau > 0$, when considering $\mathcal{R}_{u_{max}} < 1$.

\begin{figure}[H]
  \centering
  \begin{subfigure}[b]{\textwidth}
    \centering
    \includegraphics[trim={2cm 0 2cm 0},clip,width=\textwidth]{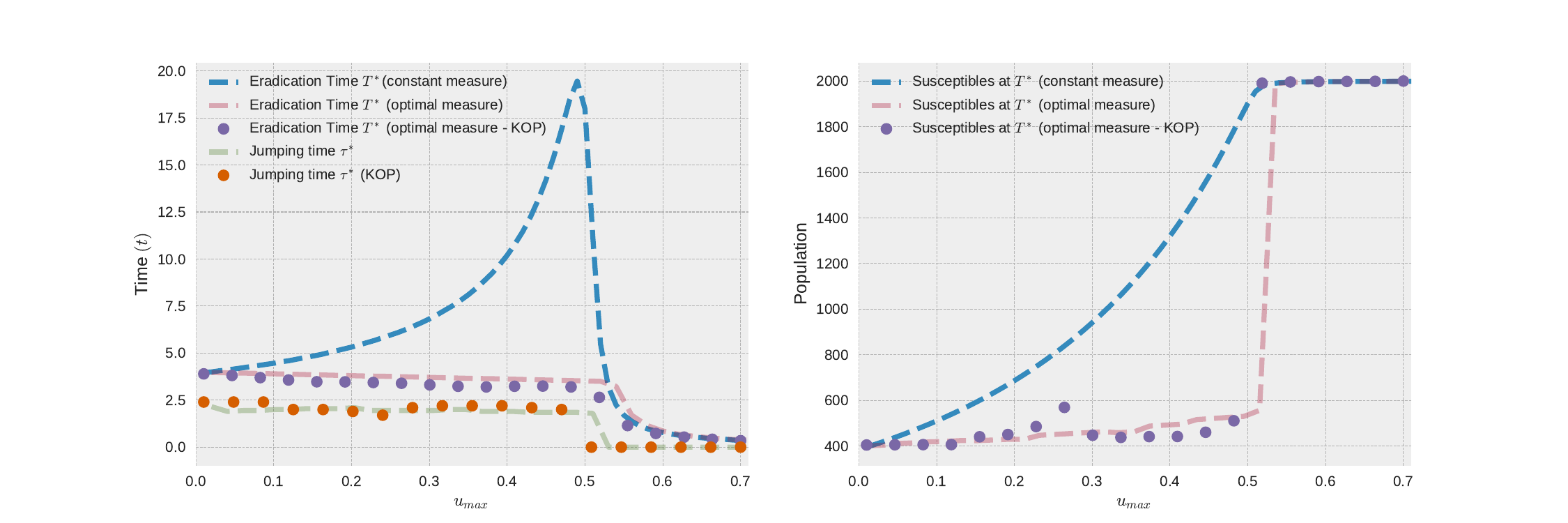}
    \caption{KOL-m.}
    \label{fig:groppisubfig-a}
  \end{subfigure}
  \\
  \begin{subfigure}[b]{\textwidth}
    \centering
    \includegraphics[trim={2cm 0 2cm 0},clip,width=\textwidth]{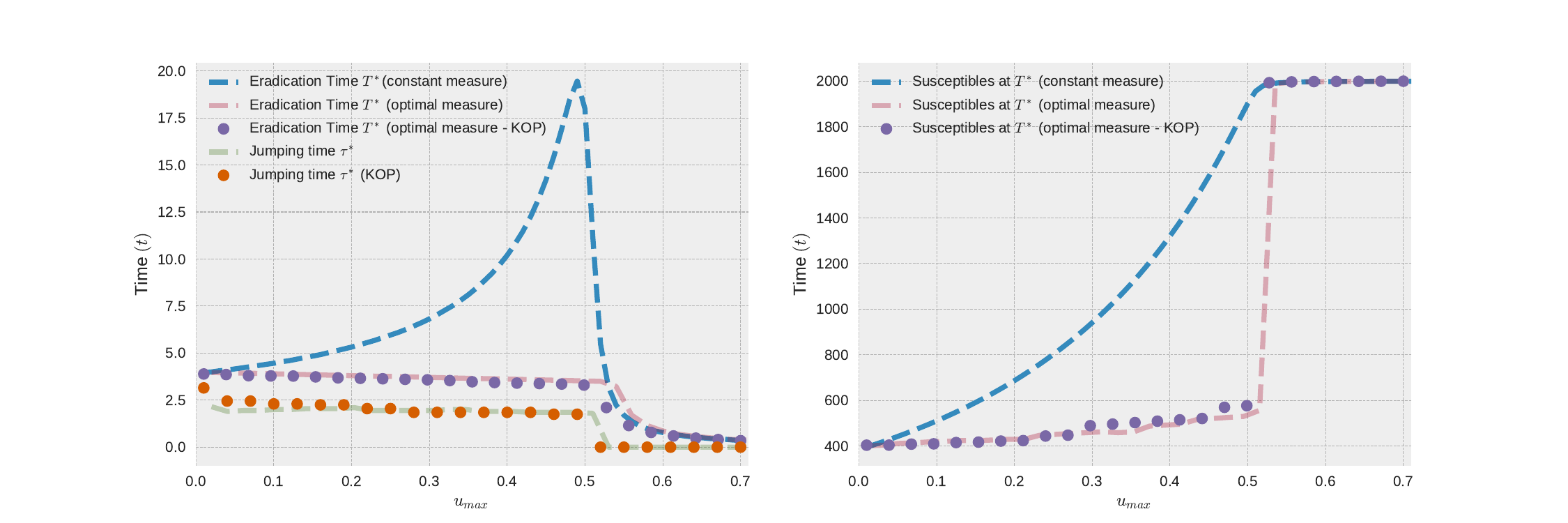}
    \caption{KOL-$\partial$.}
    \label{fig:groppisubfig-b}
  \end{subfigure}
  \caption{Solutions of the optimal control problem searching for minimum eradication time.}
  \label{fig:groppiOC}
\end{figure}

First, we numerically reproduced the results of \cite{bolzoni2017time}, by considering $\mathcal{R}_0 = 2$, $\gamma = 5$, $dt = 1$, $t^* = 100$, and $u_{max} \, \in \, [0, 0.7]$ and solving the above minimization problem constrained by the $SIR$ differential model.
In order to match the initial conditions of \cite{bolzoni2017time}, we consider $S_0 = 2000/2001$ and $I_0 = 1/2001$.
We retrieved the same behaviours of the paper of the switching time $\tau$, the eradication time $t_e$ and of the susceptibles evaluated at the eradication time (see Figure \ref{fig:groppiOC}).
We solve the optimal control problem by direct evaluations and comparison of the values of the cost functional associated to the different control strategies in $\mathcal{A}$, that are explored in an exhaustive way by sampling $\tau$ in a fine grid of step $\Delta t = 0.01$. 
This brute force strategy is not always feasible, but in the present context is motivated by the presence of the above mentioned analytical optimal solutions. 

Then, building upon the results of Section \ref{sec:KOL-numerics}, that show that our KOL approaches can be used as proxy models for a variety of epidemic differential model, we can rewrite Problem \ref{prob:groppiOC} by replacing the deterministic $SIR$ state problem with the surrogate KOL-m or KOL-$\partial$ approach.
Also in this case, the optimal solution is obtained by subsequent evaluations of control functions over the same fine grid of $\tau$.
In this way the optimal solutions obtained with our KOL approaches can be easily compared with the benchmark scenarios, thus obtaining a direct measure of the reliability of the KOL methods. 

The KOL approaches have been trained on a specific dataset constituted by 10000 Step functions of different heights, sampled uniformly in $[0, 0.8]$ (see Figure \ref{fig:u2}).

The results are collected in Figure \ref{fig:groppiOC}.
More precisely, in Figure \ref{fig:groppisubfig-a} we represent the switching times, the eradication time and the susceptibles at the eradication time for the optimal interventions of the strategy with KOL-m as surrogate state problem (point markers) together with the benchmark results (solid lines) depending on the maximum value of the intervention ($u_{max}$).
For the sake of comparison, in the same figures we represent the eradication times and the susceptibles at the eradication time for the trivial strategy implementing a constant intervention fixed at $u_{max}$ (dashed lines).
Figure \ref{fig:groppisubfig-b} represents the same benchmark quantities of \cite{bolzoni2017time} where the different quantities are retrieved by using KOL-$\partial$.
The matching of the results shows the ability of our KOL approaches to reconstruct the optimal solution without directly relying on the differential model.

\subsubsection{Optimal control for minimizing the total amount of infected }
\label{sub:ocQuad}
Here, we consider a second control problem in which the total amount of infected individuals is minimized, as detailed in the following.
\begin{ocproblem}[constrained by ODE]
  \begin{equation*}
    \min_{u \, \in \, \mathcal{U}_{ad}} \mathcal{J}_{I,u} := C_I \int_0^{t^*} I(t)^2 dt + C_u \int_0^{t^*} u(t)^2 dt,
  \end{equation*}
  subject to the state problem
  \begin{equation*}
    \begin{cases}
    \dot{S} = -\beta (1 - u) S I, \\
    \dot{I} = \beta (1-u) S I - \gamma I,\\
    \dot{R} = \gamma I, 
    \end{cases}
  \end{equation*}
  $\forall t \, \in (0,{t^*} ]$ given $[S(0), I(0), R(0)]^T = [S_0, I_0, R_0]^T$.
  The constants $C_I$, $C_u$ are positive and they balance the effects of regularization with respect to the $L^2$-infectious term.
\end{ocproblem}
This is an optimal control problem governed by the $SIR$ model where the control $u$ acts on the transmission rate and incorporates any reducing effects.
The optimal configuration obviously depends on the choice of weights $C_I$ and $C_u$.
The cost functional presents the classical Tikhonov regularization term, which has the double-purpose of increasing the convexity of the cost functional and to account for the (economic and social) burden at which higher-level NPIs come.
For what concerns the set of admissible controls, we choose the space
\begin{equation}
  \mathcal{U}_{ad} = \left \{ u(t) = \sum_{i=1}^N u_i \mathbbm{1}_{[t_{i-1}, t_i)]}(t), \, \mathrm{with} \, \{t_i\}_{i=0, N} \in (0, t^*] \; \right \},
  \label{eq:uadQuad}
\end{equation}
containing piecewise constant functions over subintervals of equal length, with the number of subintervals that is a priori chosen.
Moreover, choosing the admissible set of controls as in \eqref{eq:uadQuad} allows to immediately recast the optimal control problems as a discrete optimization problems where the unknowns correspond to the value of the piecewise constant control function in each time slab.


In our numerical exploration, we consider the scenario where $\mathcal{R}_0 = 4$, $t^* = 5$, $dt = 0.05$, initial conditions $[S_0, \, I_0, \, R_0]^T = [0.99, \, 0.01, \, 0]^T$, $u: [0,t^*] \rightarrow [0, 0.7]$ and the number of subintervals $N = \{5,10,20\}$.
Clearly, larger values of $N$ increase the complexity of the optimization problem, but at the same time they carry more flexibility in designing the optimal control strategy.
With the above choice of the parameters, we compare the solutions of the optimal control problem constrained by $SIR$ with the solutions of the following two problems, where the ODE constraint has been replaced by the surrogate operator obtained by KOL-m and KOL-$\partial$, respectively.
\begin{ocproblem}[governed by KOL-m]
 \begin{equation*}
    \min_{u \, \in \, \mathcal{U}_{ad}} \mathcal{J}_{I,u} = C_I \int_0^{t^*} I(t)^2 dt + C_u \int_0^{t^*} u(t)^2 dt,
  \end{equation*}

  where the state problem is
  \begin{equation*}
    \begin{pmatrix}
      S(t)\\
      I(t)\\
      R(t)
    \end{pmatrix} = \begin{pmatrix}
      \mathcal{\bar{G}}_{m,S}(u)(t)\\
      \mathcal{\bar{G}}_{m,I}(u)(t)\\
      \mathcal{\bar{G}}_{m,R}(u)(t)
    \end{pmatrix}  
  \end{equation*}
  $\forall t \, \in (0,T]$ given $[S(0), I(0), R(0)]^T = [S_0, I_0, R_0]^T$.
  This problem can be equivalently rewritten as an unconstrained minimization problem:
    \begin{equation*}
    \min_{u \, \in \, \mathcal{U}_{ad}} \mathcal{J}_{u} := C_I \int_0^{t^*} \mathcal{\bar{G}}_{m,I}(u)(t)^2 dt + C_u \int_0^{t^*} u(t)^2 dt.
  \end{equation*}
\end{ocproblem}
\noindent \begin{ocproblem}[governed by KOL-$\partial$]
  \begin{equation*}
    \min_{u \, \in \, \mathcal{U}_{ad}} \mathcal{J}_{I,u} = C_I \int_0^T I(t)^2 dt + C_u \int_0^T u(t)^2 dt,
  \end{equation*}
  where the state problem is
  \begin{equation*}
    \begin{pmatrix}
      \dot{S}(t)\\
      \dot{I}(t)\\
      \dot{R}(t)
    \end{pmatrix} = \begin{pmatrix}
      \mathcal{\bar{G}}_{\partial,S}(u)(t)\\
      \mathcal{\bar{G}}_{\partial,I}(u)(t)\\
      \mathcal{\bar{G}}_{\partial,R}(u)(t)
    \end{pmatrix}  
  \end{equation*}
  $\forall t \, \in (0,T]$ given $[S(0), I(0), R(0)]^T = [S_0, I_0, R_0]^T$.
  This problem can be equivalently rewritten as an unconstrained minimization problem:
    \begin{equation*}
    \min_{u \, \in \, \mathcal{U}_{ad}} \mathcal{J}_{u} := C_I \int_0^{t^*} \left ( I_0 + \int_0^{t} \mathcal{\bar{G}}_{\partial,I}(u)(\tau) d\tau \right )^2 dt + C_u \int_0^{t^*} u(t)^2 dt.
  \end{equation*}
\end{ocproblem}
Both KOL approaches have been trained with a dataset of 800 input control functions belonging to $\mathcal{U}_{ad}$ where the values of the piecewise constant control functions have been sampled from a uniform distribution in $[0, 0.8]$.
We considered different values of $C_I > 0$ and $C_u > 0$ and solved the optimization problems using the Sequential Least Square Quadratic Programming methods (SLSQP) \cite{nocedal1999numerical} as implemented in the python library Scipy \cite{virtanen2020scipy}.
This optimization scheme consists in approximating the original problem with a sequence of quadratic problems, whose objective is a second-order approximation of the Lagrangian of the original problem, and whose constraints are linearized.

\begin{figure}[H]
  \centering
  \begin{subfigure}[b]{0.3\textwidth}
    \centering
    \includegraphics[width=\textwidth]{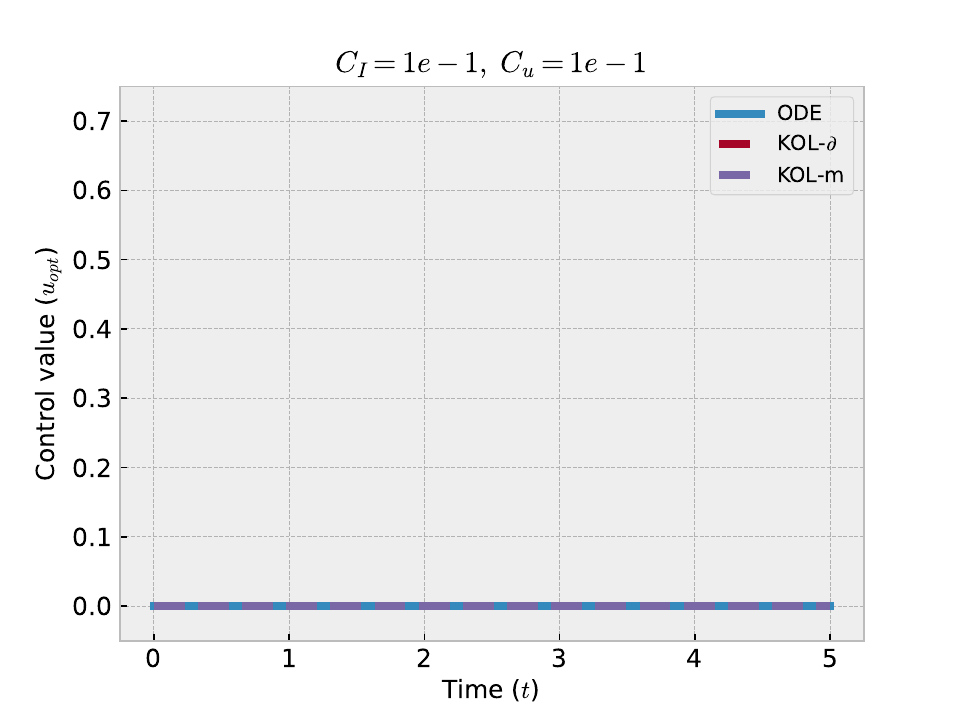}
    \caption{}
    \label{fig:cicu5subfig-a}
  \end{subfigure}
  \hfill
  \begin{subfigure}[b]{0.3\textwidth}
    \centering
    \includegraphics[width=\textwidth]{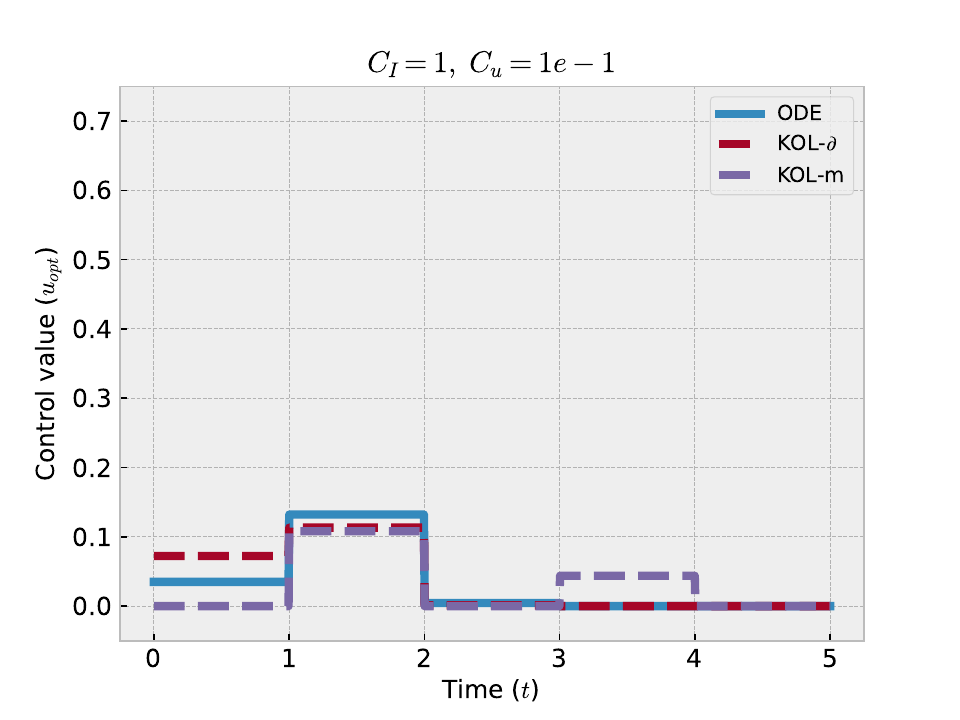}
    \caption{}
    \label{fig:cicu5subfig-b}
  \end{subfigure}
  \hfill
  \begin{subfigure}[b]{0.3\textwidth}
    \centering
    \includegraphics[width=\textwidth]{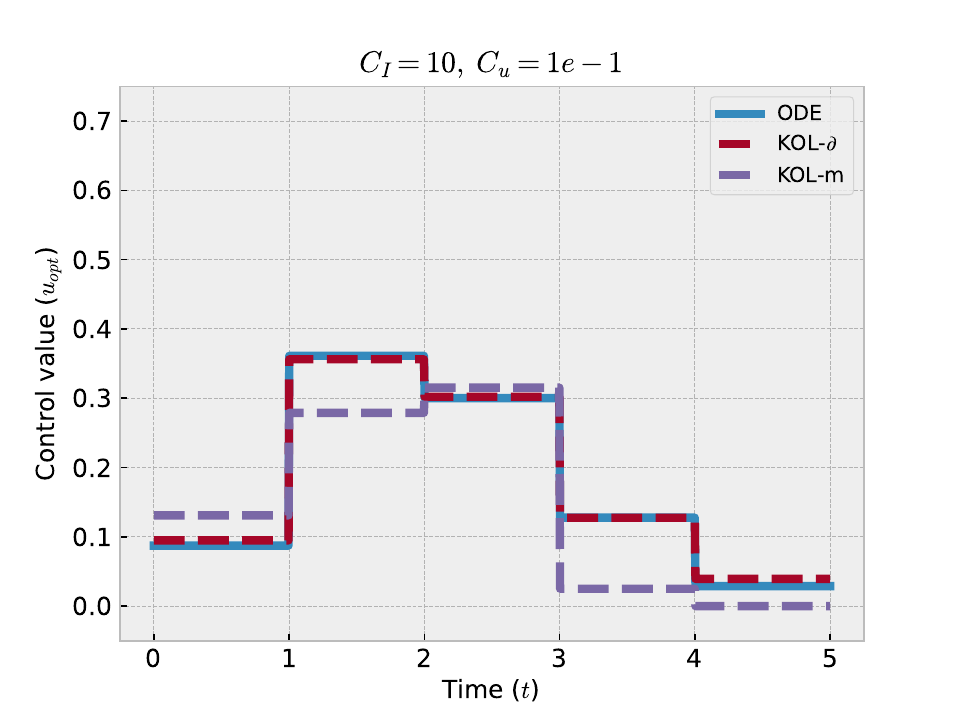}
    \caption{}
    \label{fig:cicu5subfig-c}
  \end{subfigure}
  \\
  \begin{subfigure}[b]{0.3\textwidth}
    \centering
    \includegraphics[width=\textwidth]{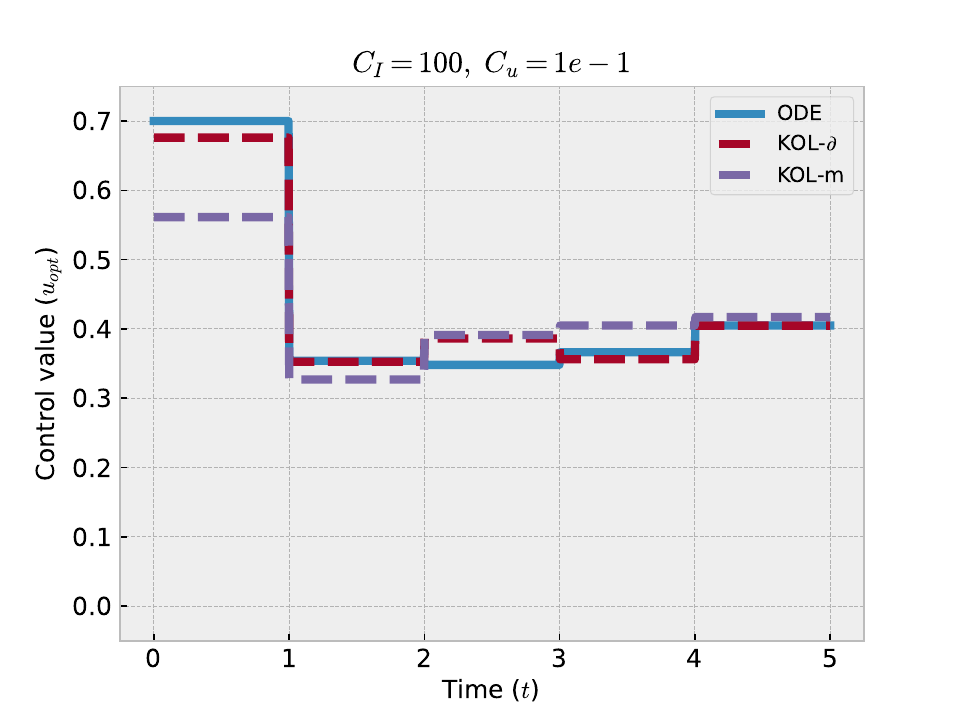}
    \caption{}
    \label{fig:cicu5subfig-d}
  \end{subfigure}
  \hfill
  \begin{subfigure}[b]{0.3\textwidth}
    \centering
    \includegraphics[width=\textwidth]{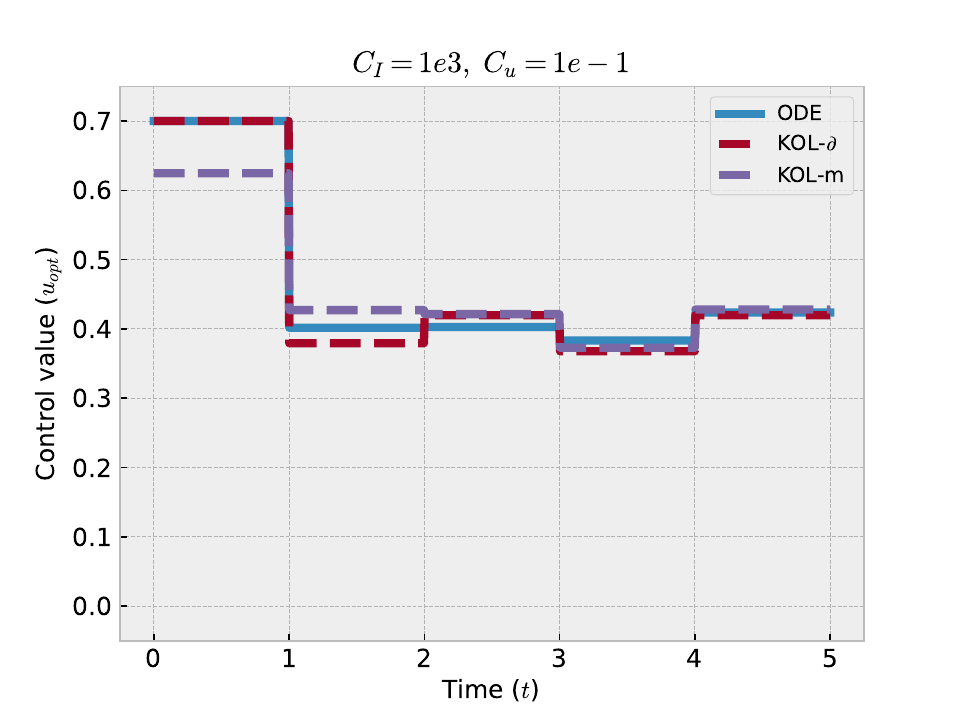}
    \caption{}
    \label{fig:cicu5subfig-e}
  \end{subfigure}
  \hfill
  \begin{subfigure}[b]{0.3\textwidth}
    \centering
    \includegraphics[width=\textwidth]{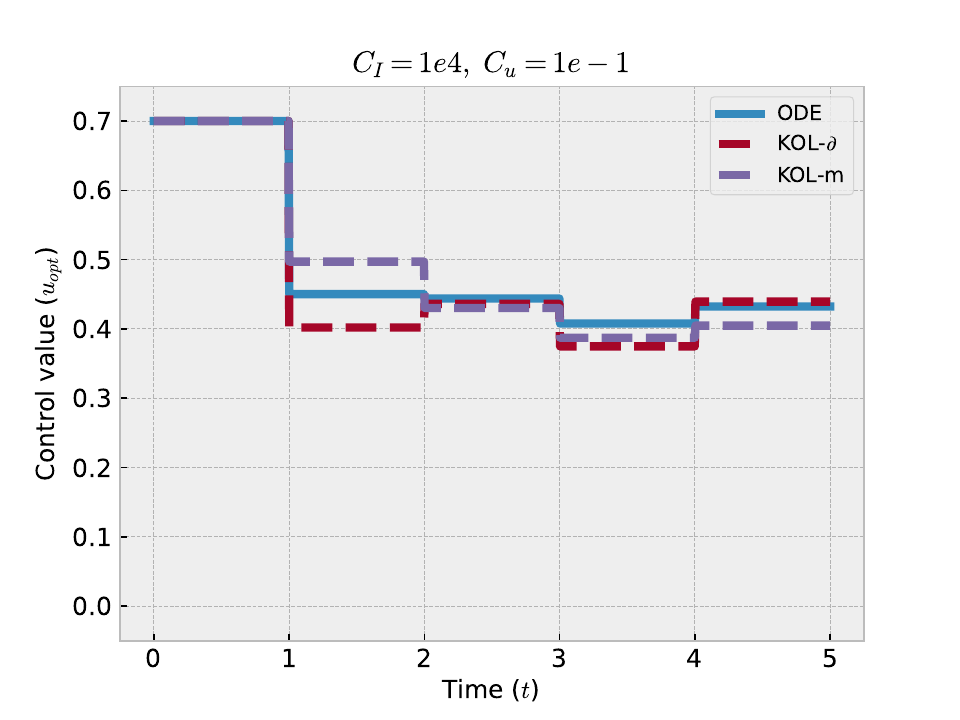}
    \caption{}
    \label{fig:cicu5subfig-f}
  \end{subfigure}
  \\
  \begin{subfigure}[b]{0.3\textwidth}
    \centering
    \includegraphics[width=\textwidth]{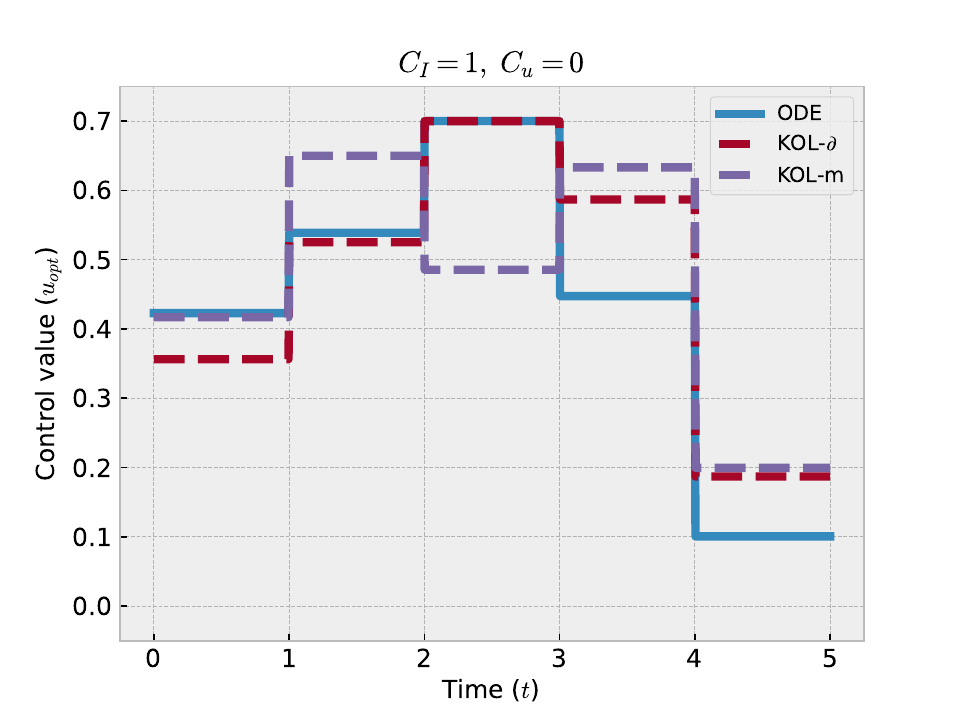}
    \caption{}
    \label{fig:cicu5subfig-g}
  \end{subfigure}
  \hfill
  \begin{subfigure}[b]{0.3\textwidth}
    \centering
    \includegraphics[width=\textwidth]{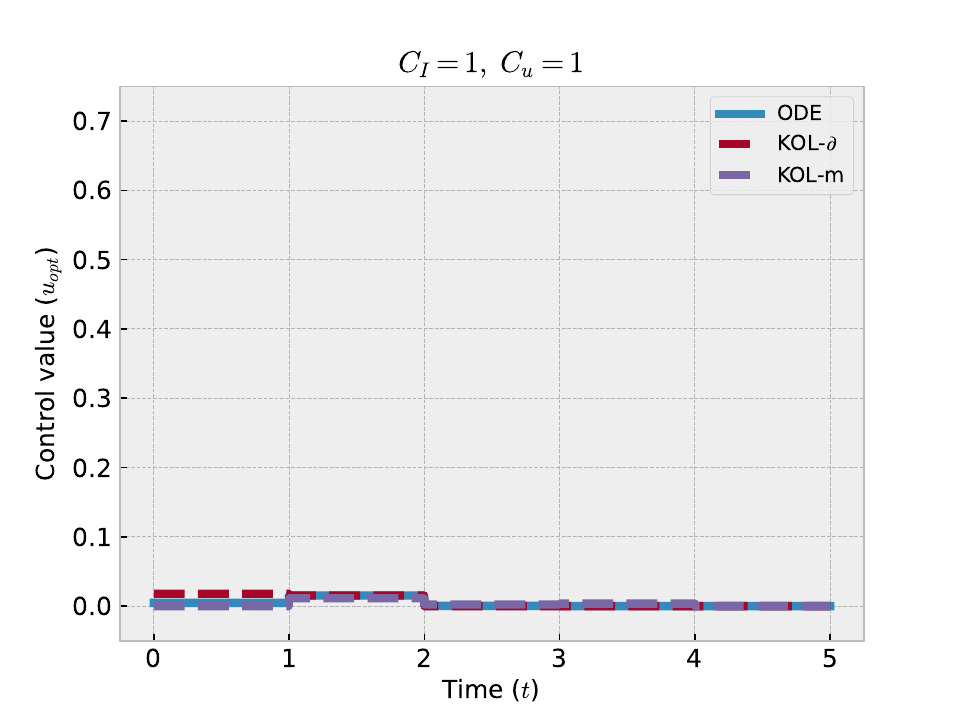}
    \caption{}
    \label{fig:cicu5subfig-h}
  \end{subfigure}
  \hfill
  \begin{subfigure}[b]{0.3\textwidth}
    \centering
    \includegraphics[width=\textwidth]{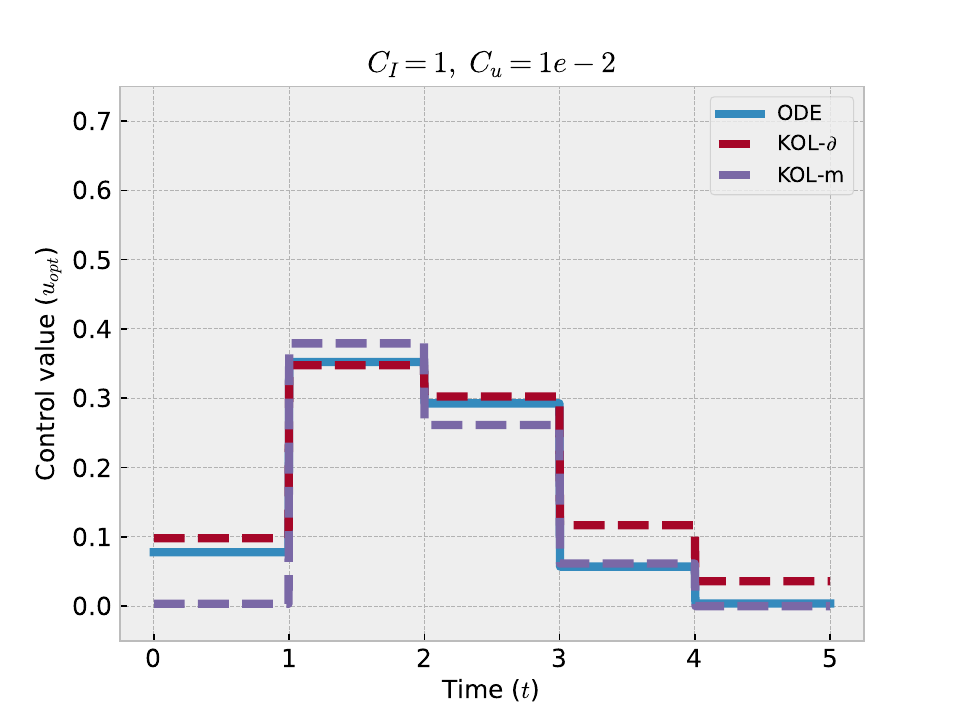}
    \caption{}
    \label{fig:cicu5subfig-i}
  \end{subfigure}
  \\
  \begin{subfigure}[b]{0.3\textwidth}
    \centering
    \includegraphics[width=\textwidth]{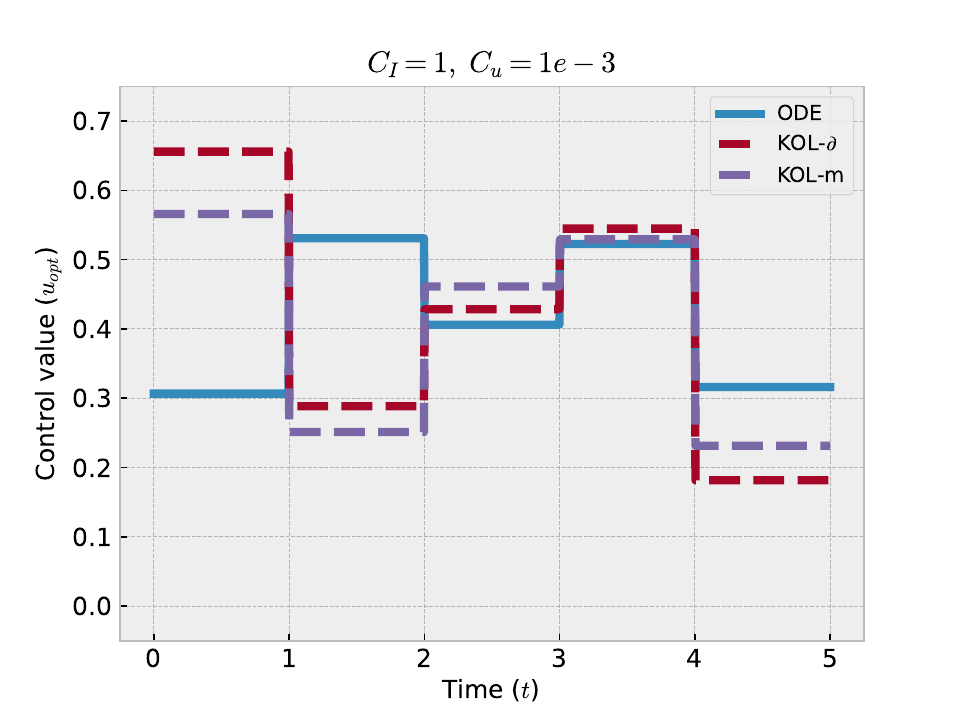}
    \caption{}
    \label{fig:cicu5subfig-l}
  \end{subfigure}
  \hfill
  \begin{subfigure}[b]{0.3\textwidth}
    \centering
    \includegraphics[width=\textwidth]{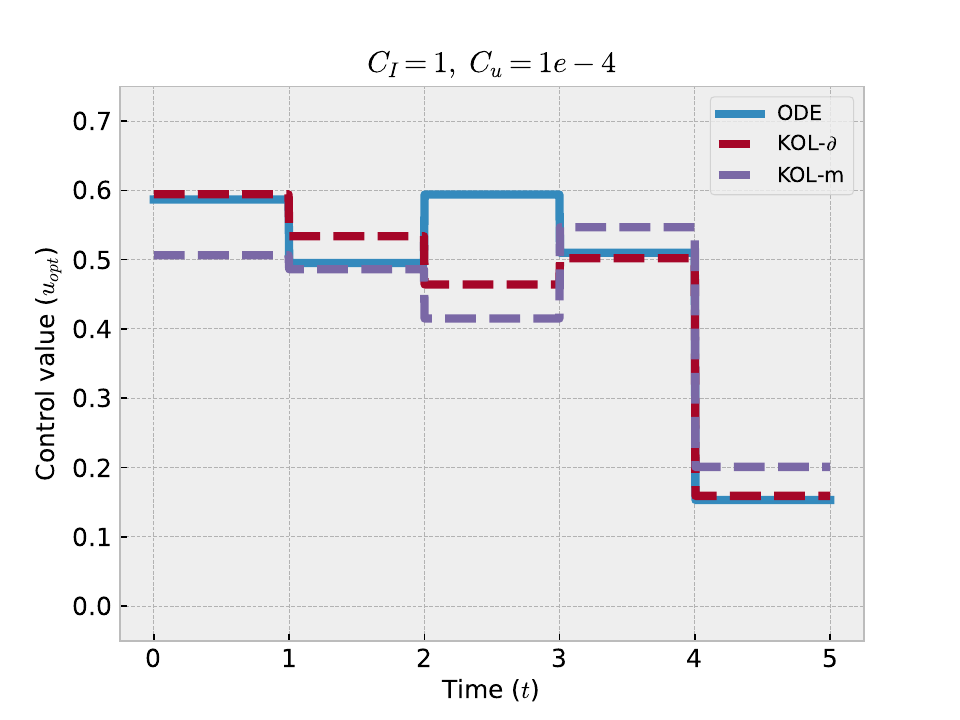}
    \caption{}
    \label{fig:cicu5subfig-m}
  \end{subfigure}
  \hfill
  \begin{subfigure}[b]{0.3\textwidth}
    \centering
    \includegraphics[width=\textwidth]{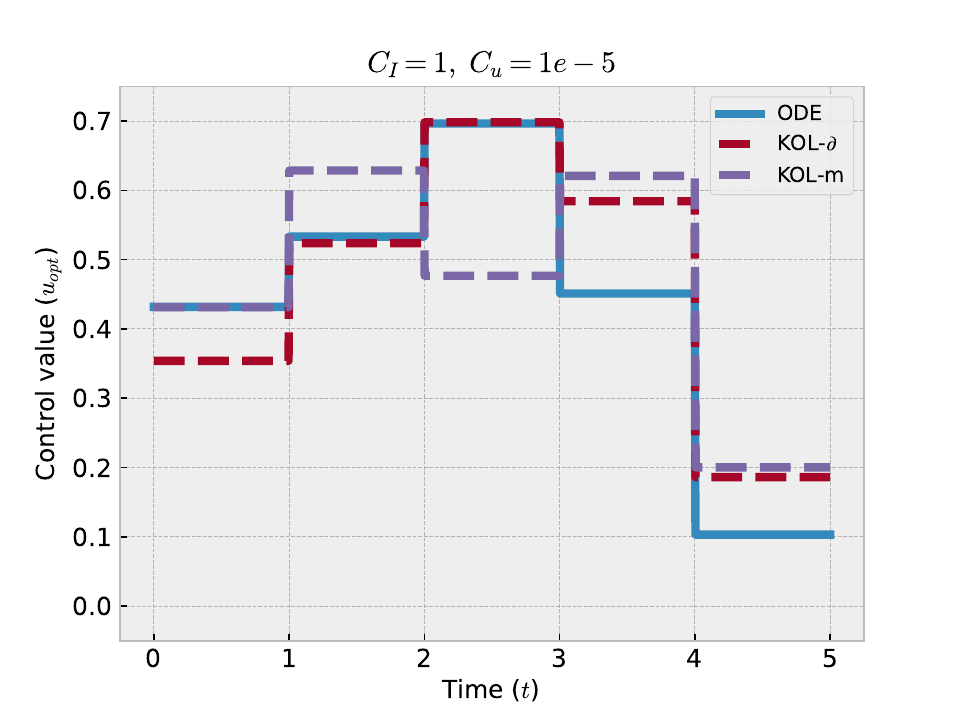}
    \caption{}
    \label{fig:cicu5subfig-n}
  \end{subfigure}
  \\
  \begin{subfigure}[b]{0.3\textwidth}
    \centering
    \includegraphics[width=\textwidth]{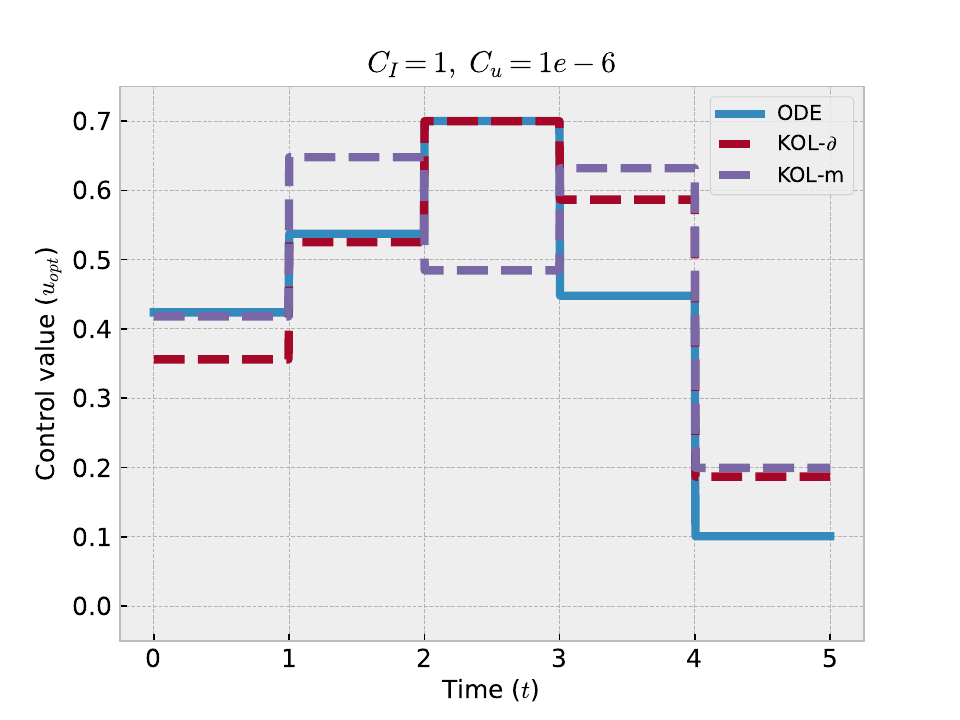}
    \caption{}
    \label{fig:cicu5subfig-o}
  \end{subfigure}
  \caption{Optimal controls for the three optimal control problems fixing $N=5$.}
  \label{fig:cicu5}
\end{figure}

\begin{figure}[H]
  \centering
  \begin{subfigure}[b]{0.45\textwidth}
    \includegraphics[width=\textwidth]{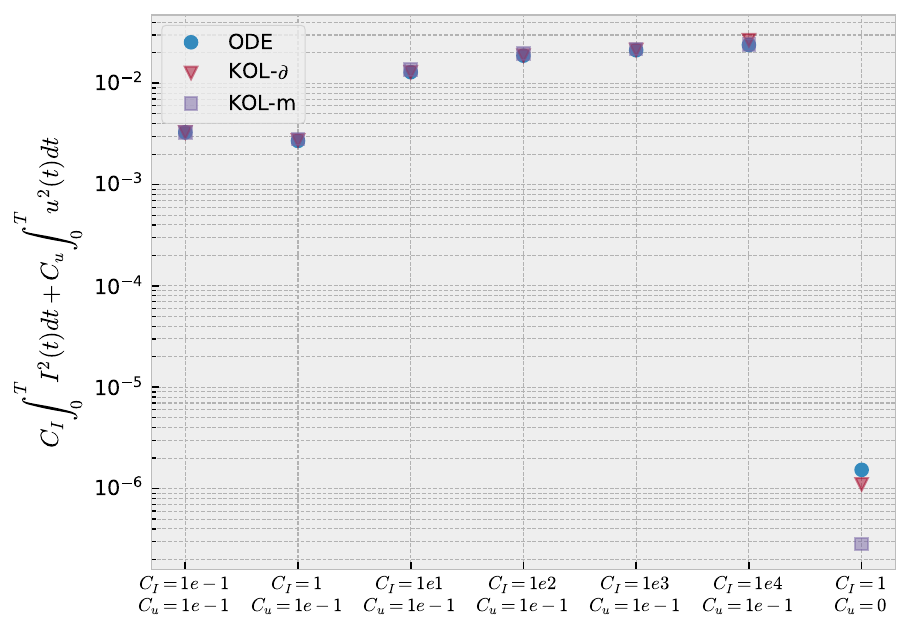}
    \caption{}
    \label{fig:subfigs_1}
  \end{subfigure}
  \hfill
  \begin{subfigure}[b]{0.45\textwidth}
    \includegraphics[width=\textwidth]{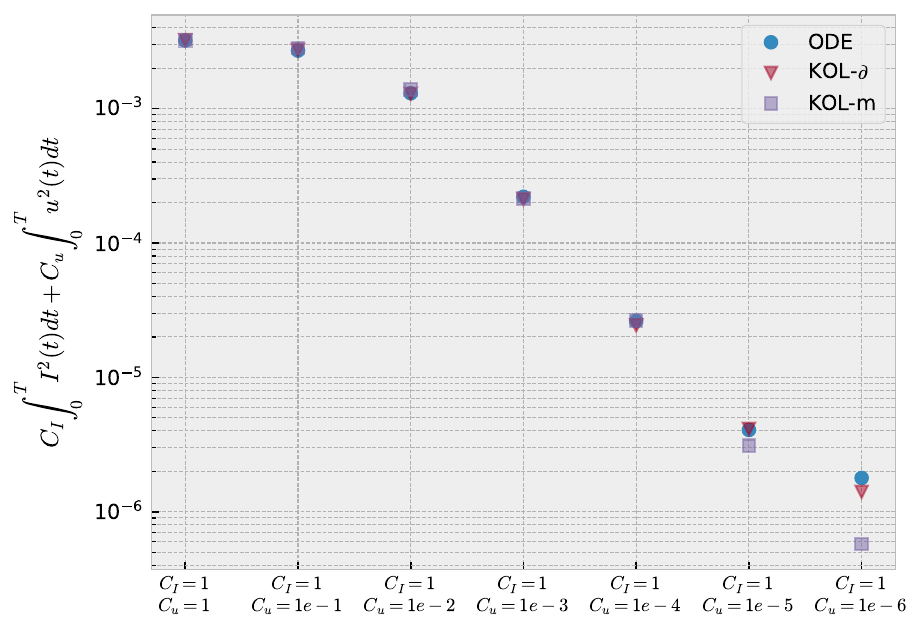}
    \caption{}
    \label{fig:subfigs_2}
  \end{subfigure}
  \caption{Cost functionals at the optimal control for different $C_I$ and $C_u$ ($N=5$). (a) We consider different orders of magnitude for $C_I$, keeping $C_u = 1e-1$. (b) We consider different orders of magnitude for $C_u$, keeping $C_I = 1$. }
  \label{fig:errorsOCquad5}
\end{figure}

In the case of $N=5$ phases, Figure \ref{fig:cicu5} shows 13 different scenarios of the optimal trajectories reconstructed by solving the problem under $SIR$ (solid line), KOL-m (purple dashed line) and KOL-$\partial$ (red dashed line) constraint, respectively.
At first glance, we see that the optimal controls obtained with the three approaches do not exactly coincide, especially when $C_u \ll C_I $, \textit{i.e.} when the problem is extremely non-convex.
However, comparing the cost functional values associated to the three optimal controls shed a completely different light on the reliability of our approach.
Indeed, solving the deterministic $SIR$ problem for each optimal control strategy, and plotting the corresponding cost functional values in Figure \ref{fig:errorsOCquad5}, we can appreciate that the cost values practically coincide, thus showing the efficacy of our KOL-based approach in solving the optimal control problem. 
It is worth remarking that in some cases, typically the less convex ones (i.e. $C_u$ much larger than $C_I $), the KOL approaches succeed in finding cost values lower than the ones obtained using the $SIR$ constraint.

\begin{figure}[H]
  \centering
  \begin{subfigure}[b]{0.3\textwidth}
    \centering
    \includegraphics[width=\textwidth]{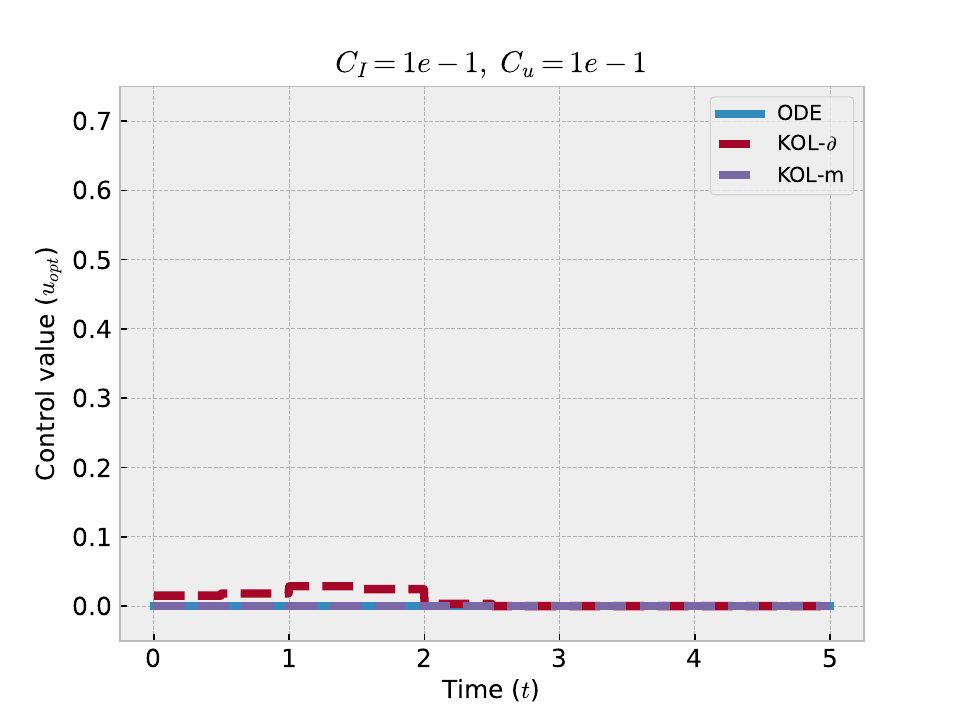}
    \caption{}
    \label{fig:cicu10subfig-a}
  \end{subfigure}
  \hfill
  \begin{subfigure}[b]{0.3\textwidth}
    \centering
    \includegraphics[width=\textwidth]{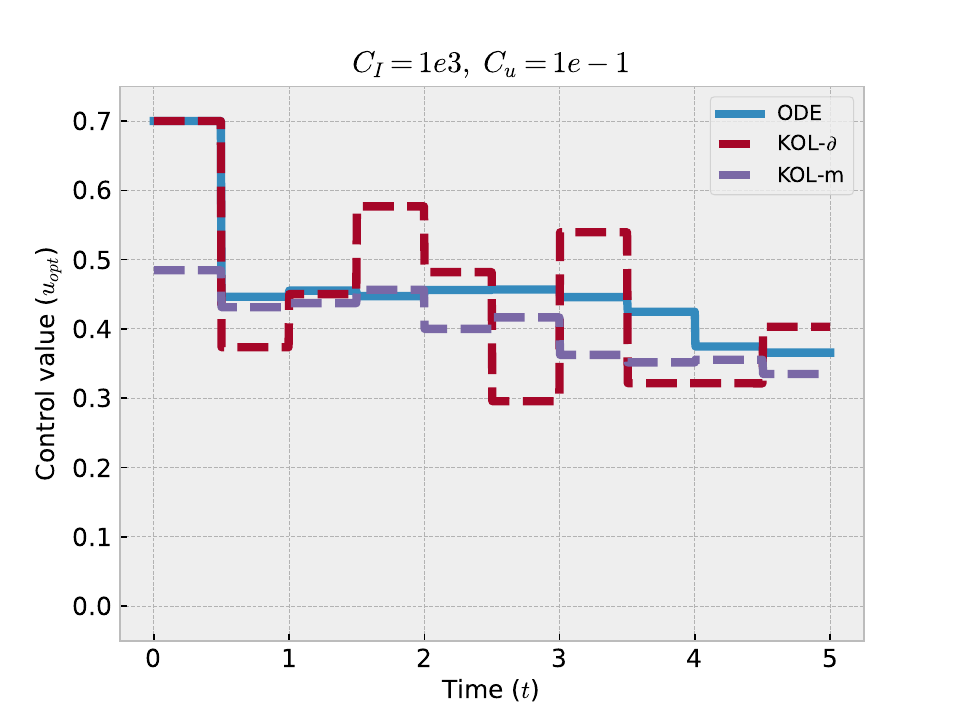}
    \caption{}
    \label{fig:cicu10subfig-b}
  \end{subfigure}
  \hfill
  \begin{subfigure}[b]{0.3\textwidth}
    \centering
    \includegraphics[width=\textwidth]{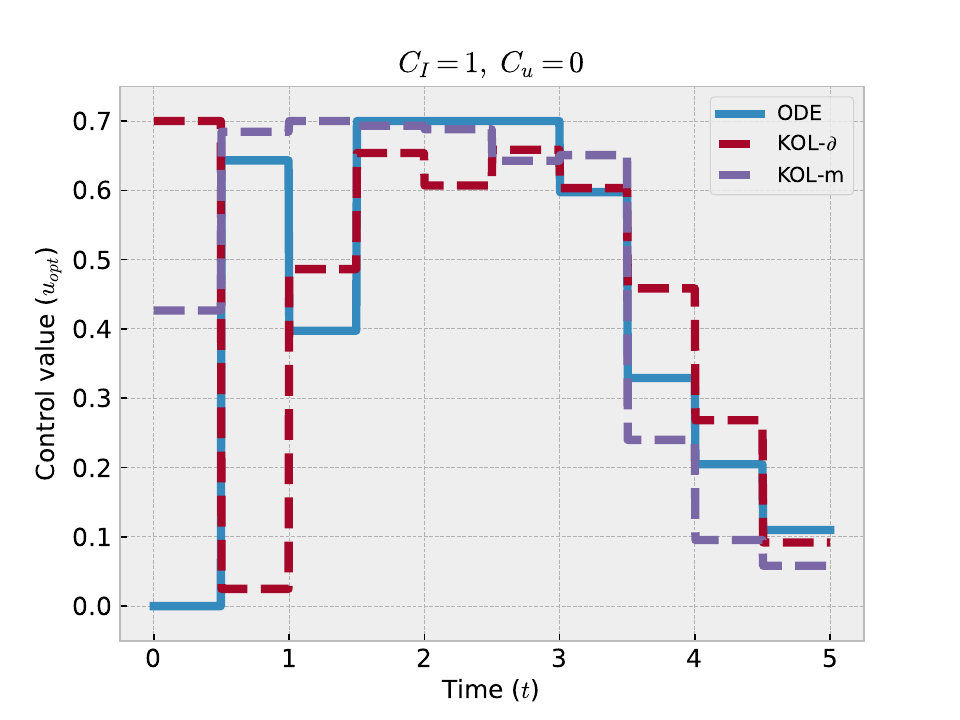}
    \caption{}
    \label{fig:cicu10subfig-c}
  \end{subfigure}
  \\
  \begin{subfigure}[b]{0.3\textwidth}
    \centering
    \includegraphics[width=\textwidth]{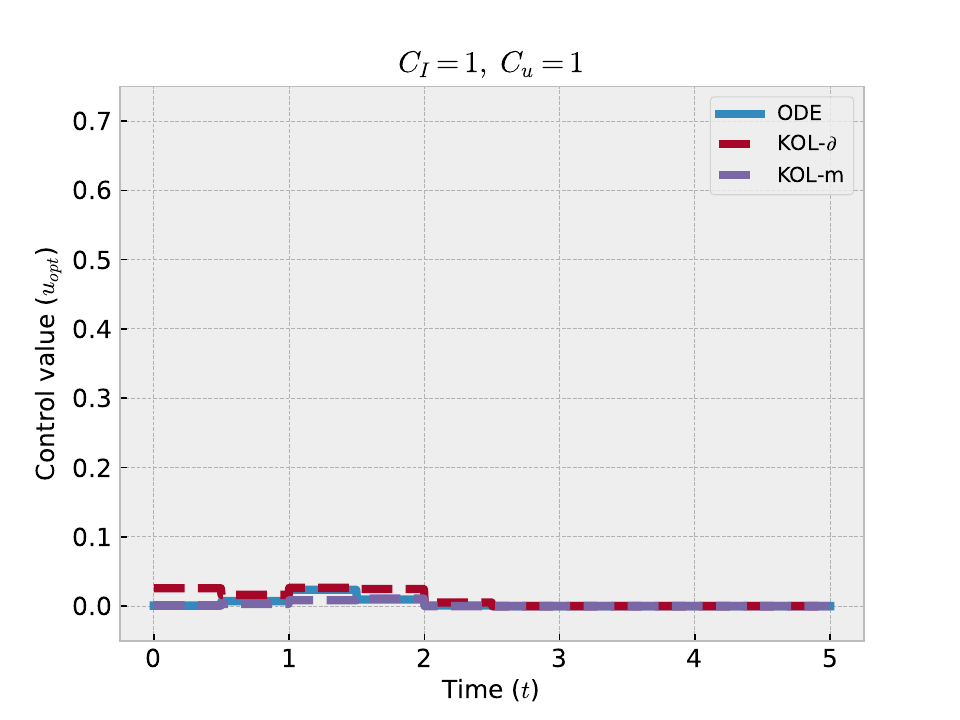}
    \caption{}
    \label{fig:cicu10subfig-d}
  \end{subfigure}
  \hfill
  \begin{subfigure}[b]{0.3\textwidth}
    \centering
    \includegraphics[width=\textwidth]{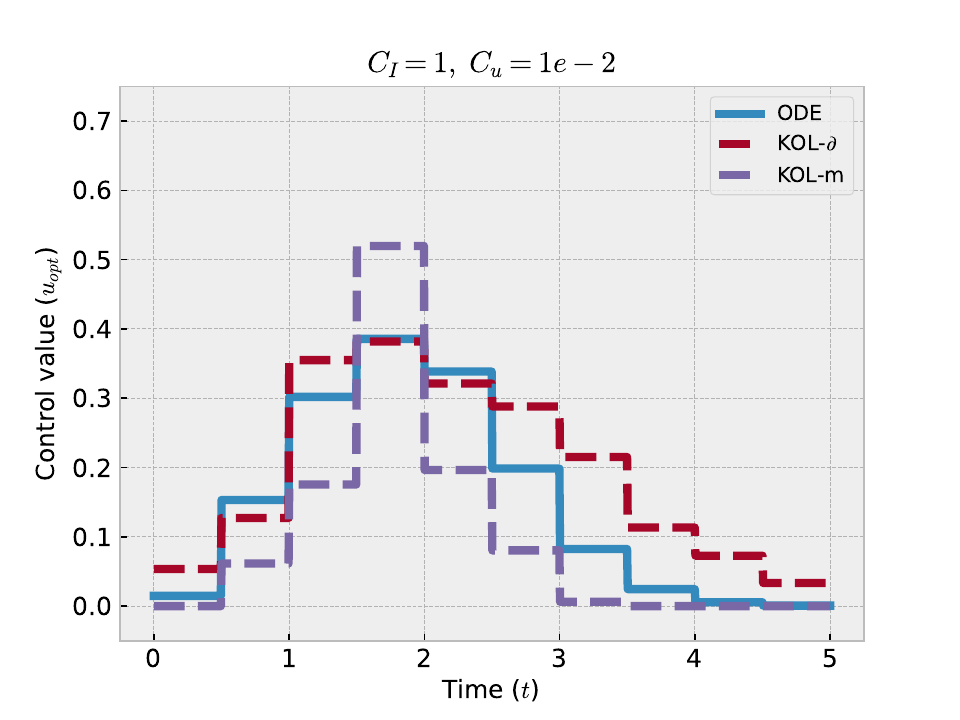}
    \caption{}
    \label{fig:cicu10subfig-e}
  \end{subfigure}
  \hfill
  \begin{subfigure}[b]{0.3\textwidth}
    \centering
    \includegraphics[width=\textwidth]{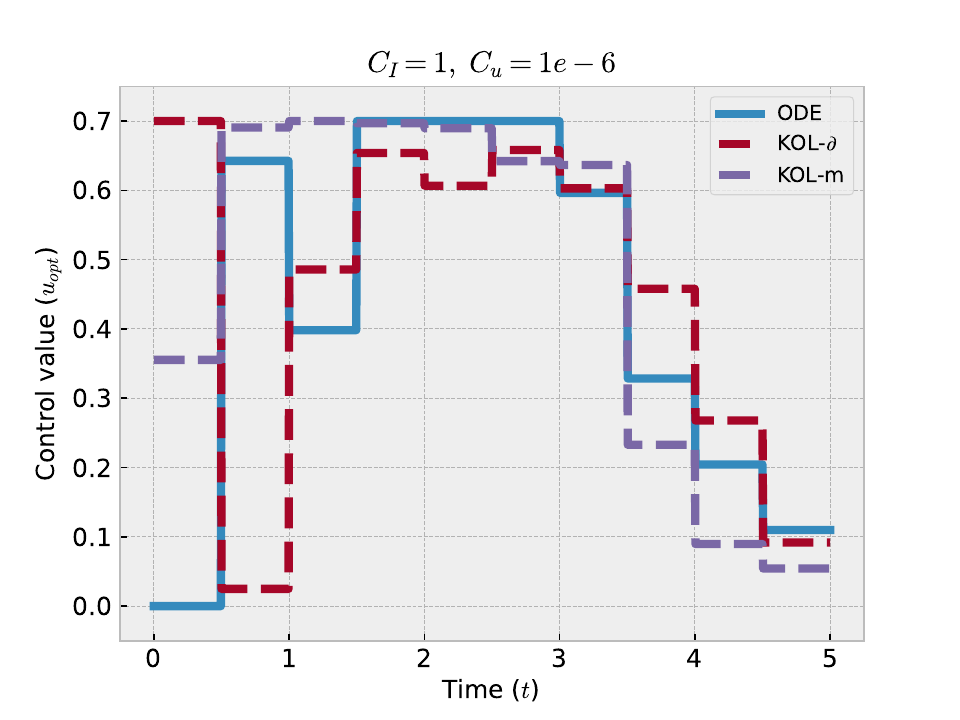}
    \caption{}
    \label{fig:cicu10subfig-f}
  \end{subfigure}  
  \caption{Optimal controls for the three optimal control problems fixing $N=10$.}
  \label{fig:cicu10subfigures}
\end{figure}

\begin{figure}[H]
  \centering
  \begin{subfigure}[b]{0.3\textwidth}
    \centering
    \includegraphics[width=\textwidth]{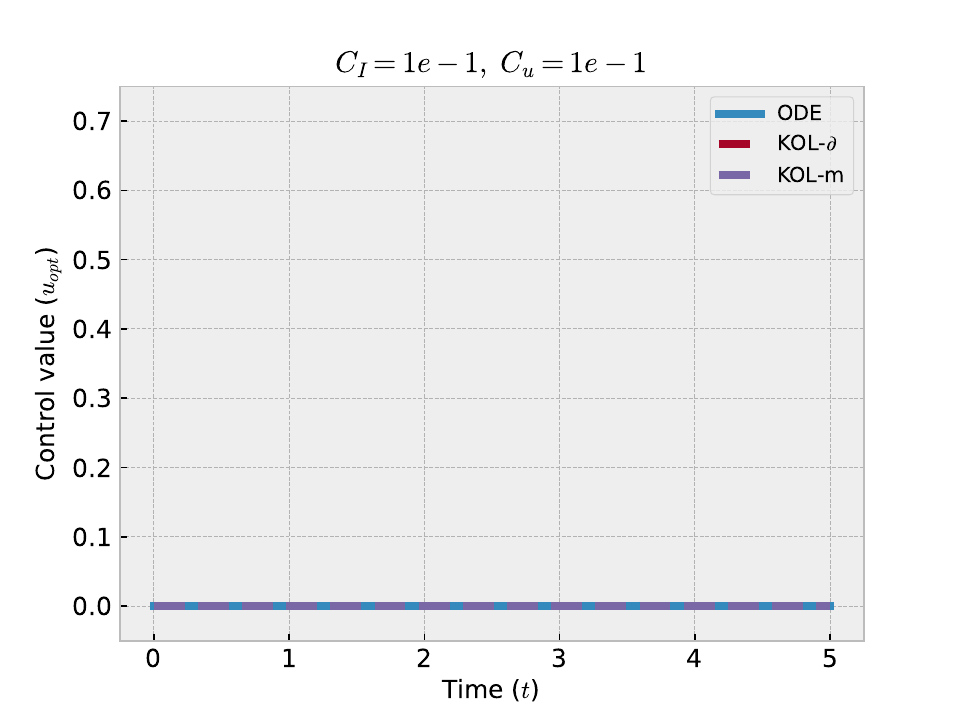}
    \caption{}
    \label{fig:cicu20subfig-a}
  \end{subfigure}
  \hfill
  \begin{subfigure}[b]{0.3\textwidth}
    \centering
    \includegraphics[width=\textwidth]{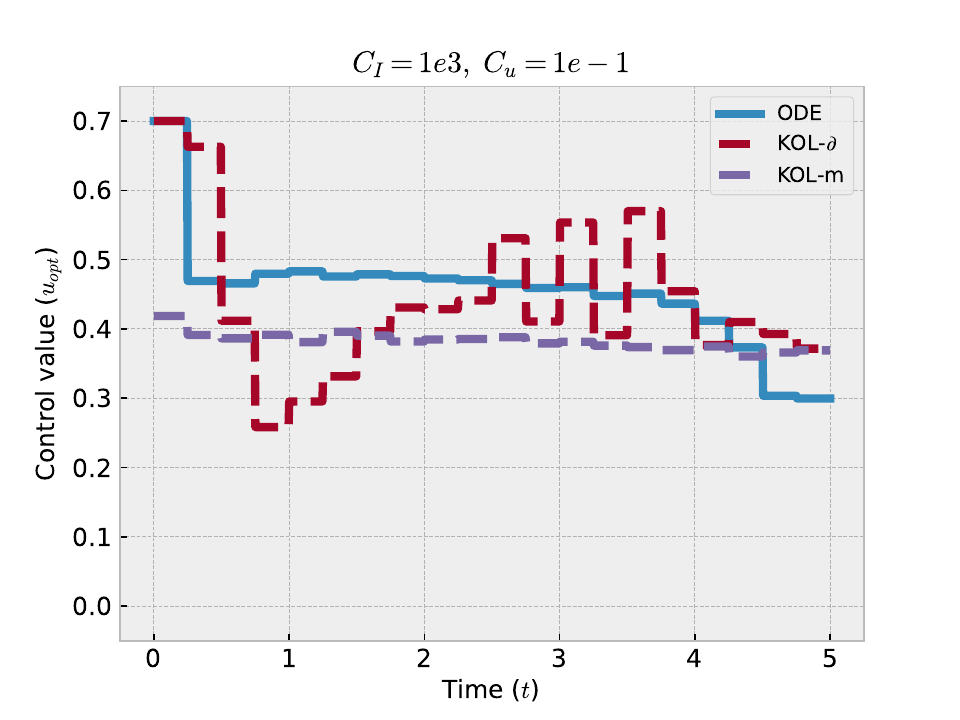}
    \caption{}
    \label{fig:cicu20subfig-b}
  \end{subfigure}
  \hfill
  \begin{subfigure}[b]{0.3\textwidth}
    \centering
    \includegraphics[width=\textwidth]{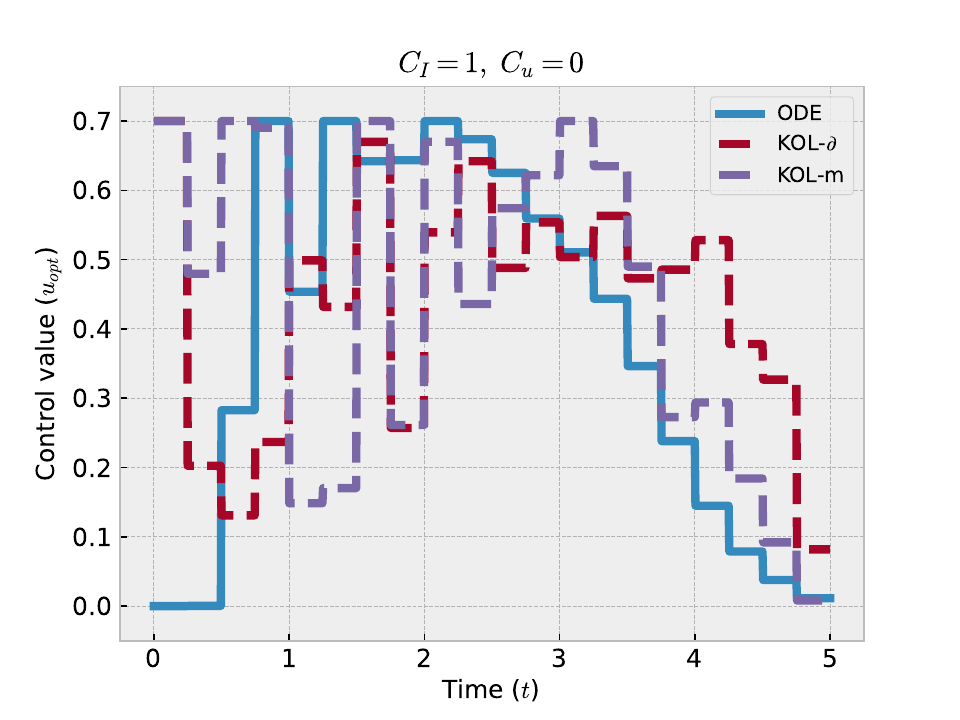}
    \caption{}
    \label{fig:cicu20subfig-c}
  \end{subfigure}
  \\
    \begin{subfigure}[b]{0.3\textwidth}
    \centering
    \includegraphics[width=\textwidth]{CI1e_1Cu1e_1_20.pdf}
    \caption{}
    \label{fig:cicu20subfig-d}
  \end{subfigure}
  \hfill
  \begin{subfigure}[b]{0.3\textwidth}
    \centering
    \includegraphics[width=\textwidth]{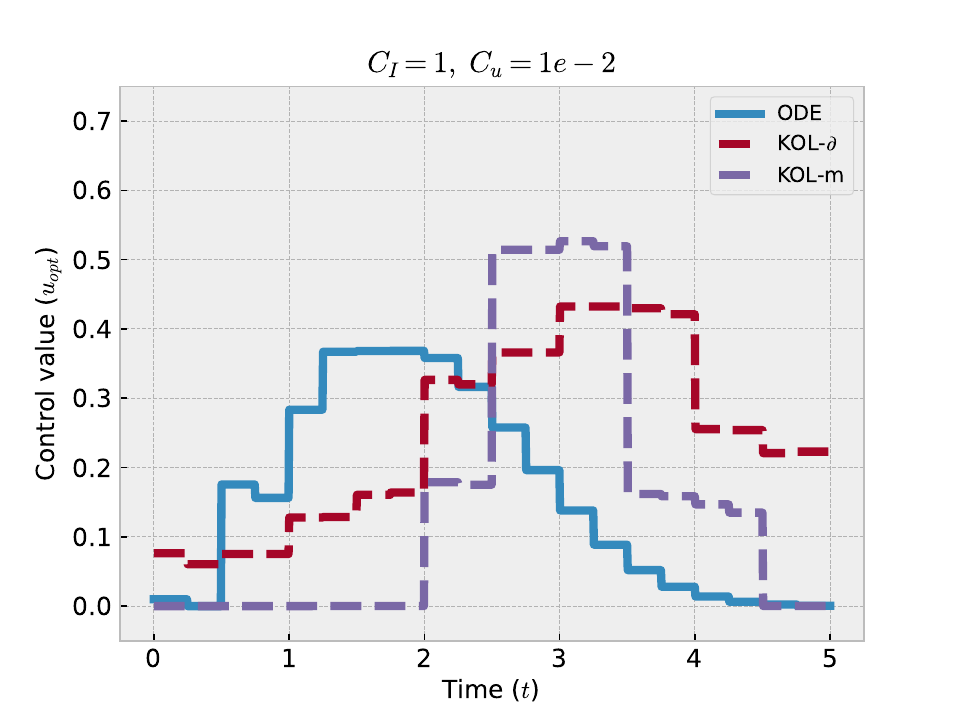}
    \caption{}
    \label{fig:cicu20subfig-e}
  \end{subfigure}
  \hfill
  \begin{subfigure}[b]{0.3\textwidth}
    \centering
    \includegraphics[width=\textwidth]{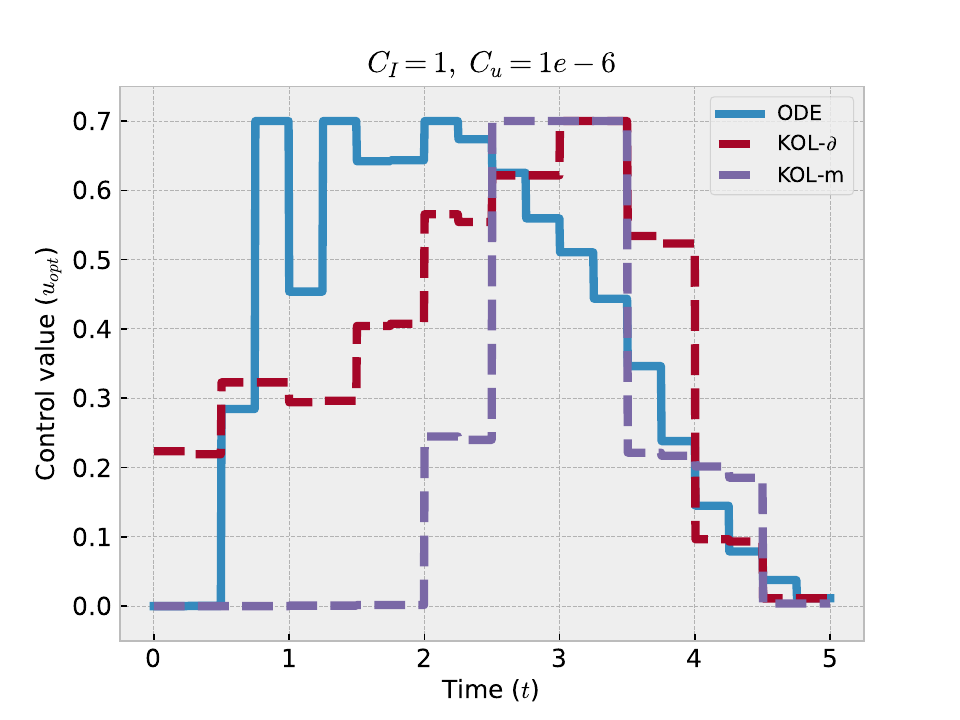}
    \caption{}
    \label{fig:cicu20subfig-f}
  \end{subfigure}
  \caption{Optimal controls for the three optimal control problems fixing $N=20$.}
  \label{fig:cicu20subfigures}
\end{figure}

\begin{figure}[H]
  \centering
  \begin{subfigure}[b]{0.45\textwidth}
    \includegraphics[width=\textwidth]{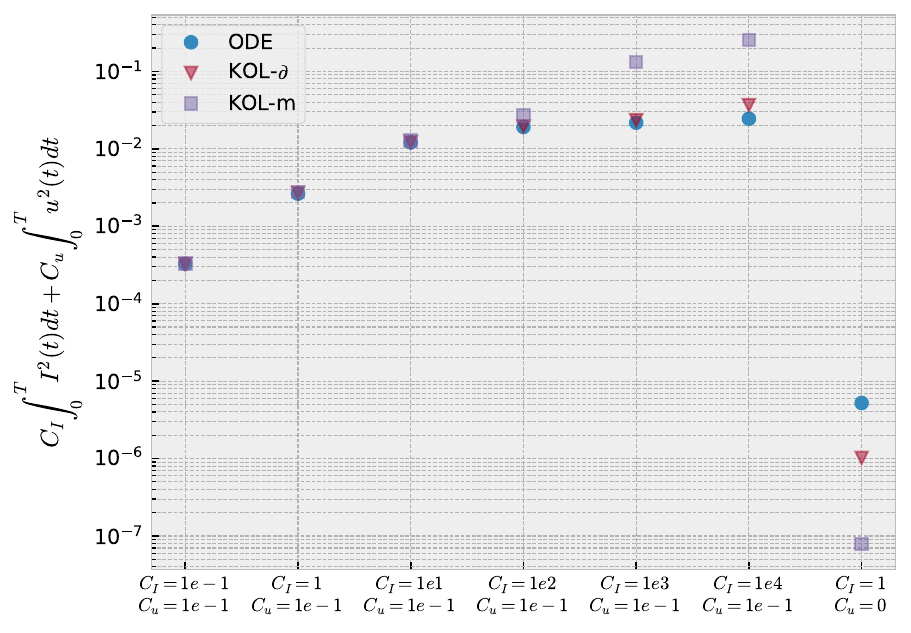}
    \caption{}
    \label{fig:subfigs_2a}
  \end{subfigure}
  \hfill
  \begin{subfigure}[b]{0.45\textwidth}
    \includegraphics[width=\textwidth]{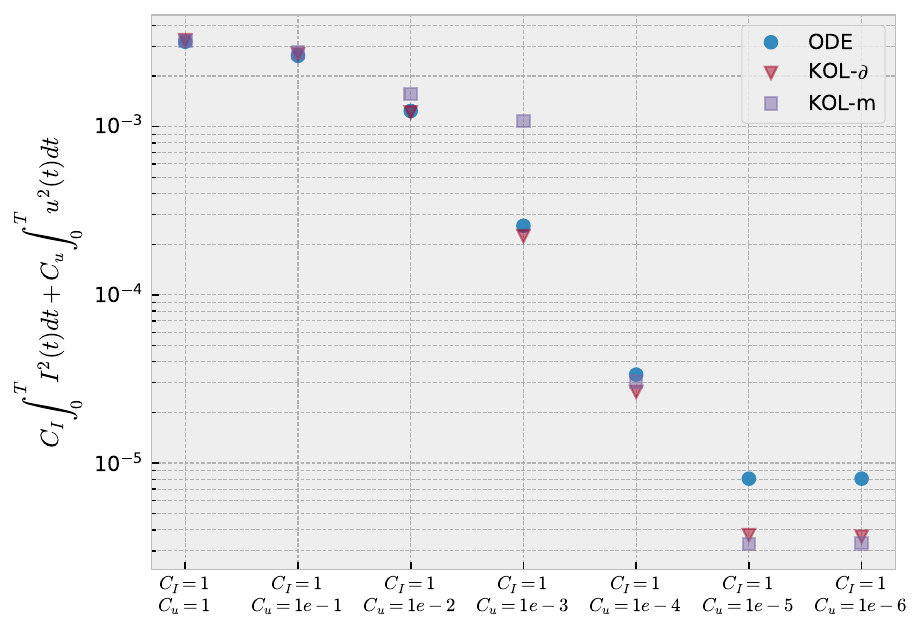}
    \caption{}
    \label{fig:subfigs_2b}
  \end{subfigure}
  \caption{Cost functionals at the optimal control for different $C_I$ and $C_u$ ($N=10$). (a) We consider different orders of magnitude for $C_I$, keeping $C_u = 1e-1$. (b) We consider different orders of magnitude for $C_u$, keeping $C_I = 1$.}
  \label{fig:errorsOCquad10}
\end{figure}

\begin{figure}[H]
  \centering
  \begin{subfigure}[b]{0.45\textwidth}\includegraphics[width=\textwidth]{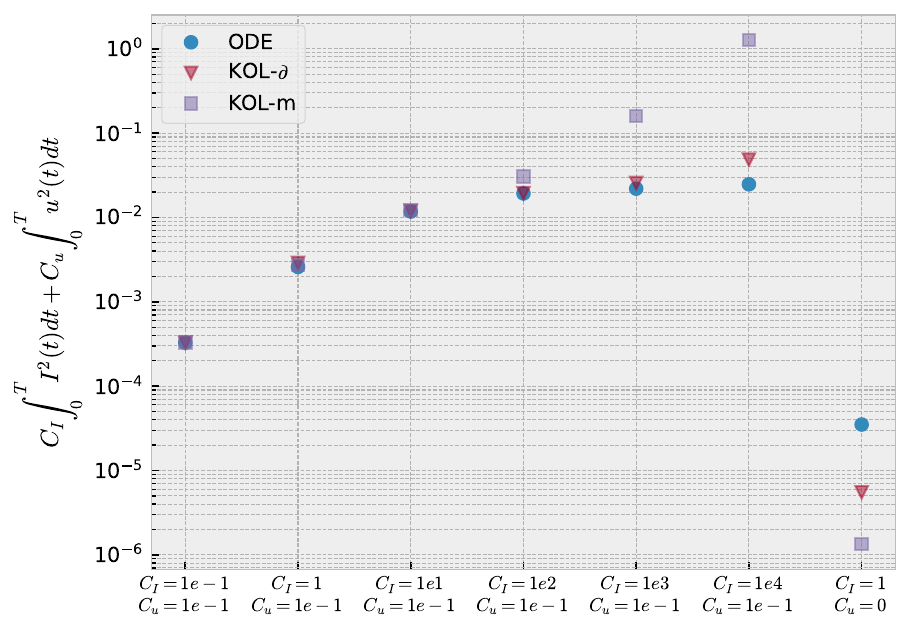}
  \caption{}
  \label{fig:subfig3a}
  \end{subfigure}
  \hfill
  \begin{subfigure}[b]{0.45\textwidth}\includegraphics[width=\textwidth]{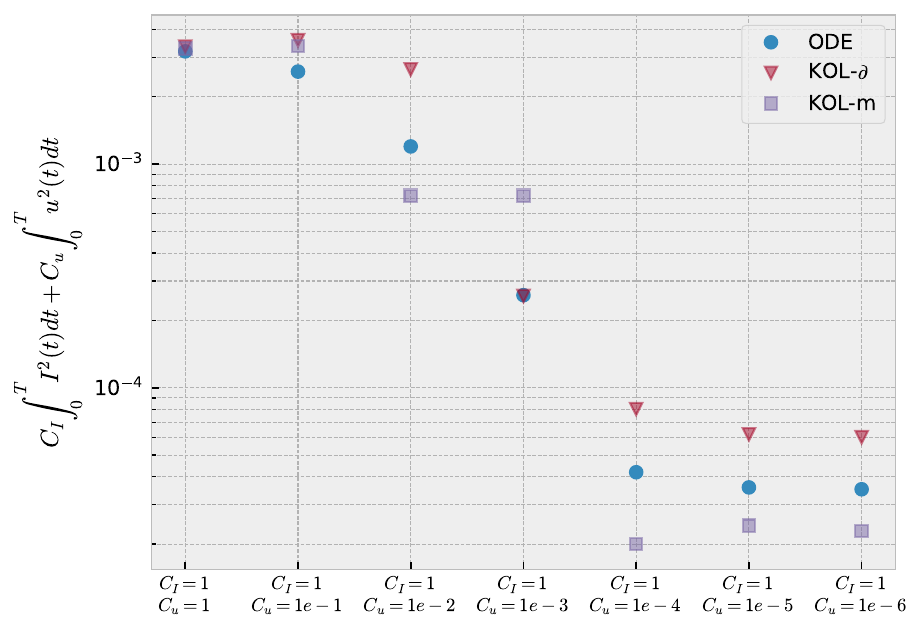}
  \caption{}
  \label{fig:subfig3b}
  \end{subfigure}
  \caption{Cost functionals at the optimal control for different $C_I$ and $C_u$ ($N=20$). (a) We consider different orders of magnitude for $C_I$, keeping $C_u = 1e-1$. (b) We consider different orders of magnitude for $C_u$, keeping $C_I = 1$.}
  \label{fig:errorsOCquad20}
\end{figure}

Finally, we compare the three approaches for increasing complexity of the optimization problem, namely for $N=10,20$.
The results are collected in Figures
\ref{fig:cicu10subfigures}-\ref{fig:cicu20subfigures},
where the optimal control functions are reported, and in Figures \ref{fig:errorsOCquad10}-\ref{fig:errorsOCquad20}, where the corresponding cost values are compared. 
Also in these cases, we observe that the cost values achieved by the two KOL approaches are comparable, if not lower, than those obtained employing the $SIR$ constraint.
However, with a growing amount of phases the three approaches produce different optimal trajectories due to the increased non-convexity of the cost functional and, in cascade, of the optimal control problem itself, which exhibit many local minima.

This latest set of results further proves the reliability and the effectiveness of KOL approaches for the most common optimization tasks.
\section*{Acknowledgements}
The authors thank Dr. F. Regazzoni for the support provided with \NNML algorithm and for the insightful discussions on the topic.
MV and GZ are members of INdAM-GNCS.
The present reasearch is part of the activities of "Dipartimento di Eccellenza 2023-2027", MUR, Italy.

\section*{Data Availability}
The computational framework is available on GitHub together with the numerical results of this work: \texttt{https://github.com/giovanniziarelli/KOL}.

\clearpage
\bibliographystyle{vancouver}
\bibliography{main}
\newpage
\appendix
\renewcommand{\thesection}{\Alph{section}}
\section{Summary of Neural Tangent Kernel}
\label{app:ntk}
  We introduce the Neural Tangent Kernel in the specific context of neural networks minimizing the Mean Square Error in the input dataset.
  Consider a set of input data (training set)
  \begin{equation}
    \mathcal{D} = \{ (x_1,y_1), (x_2,y_2) \hdots (x_n, y_n)\},
  \end{equation}
  where $\{x_i\}_i \in \mathbb{R}^d$ represents the sequence of input data, and $\{ y_i\}_i \in \mathbb{R}$ the respective labels.
  Let $f(x; \theta)$ be the neural network regressor, whose weights $\theta$ need to be trained solving the minimization problem
  \begin{equation}
    \min_{\theta \in \mathbb{R}^{n_{\theta}}} L(\theta) = \sum_{i=1}^n \dfrac{1}{2} (f(x_i; \theta) - y_i)^2 = \sum_{i=1}^n l(f(x_i; \theta), y_i).
  \end{equation}
  We solve the optimization problem via gradient flow, which is the continuous counterpart of the usual full-batch gradient descent method used in machine learning,
  \begin{equation}
    \partial_\tau \theta(\tau) = - \partial_{\theta} L (\theta(\tau)) = -\sum_{i=1}^n \partial_{\theta} f(x_i; \theta(\tau)) (f(x_i; \theta(\tau)) - y_i),
  \end{equation}
  where $\tau$ is a fictious time variable accounting for the iteration progress.
  Therefore, 
  \begin{equation}
  \begin{split}
    \partial_\tau f(x_j; \theta(\tau)) &= \partial_{\theta} f(x_j; \theta(\tau)) \partial_\tau \theta(\tau) \\
    &= -\sum_{i=1}^n \partial_{\theta} f(x_j; \theta(\tau))^T \partial_{\theta} f(x_i; \theta(\tau)) (f(x_i; \theta(\tau)) - y_i), 
  \end{split}
  \end{equation}
  and we define the Neural Tangent Kernel $K: \mathbb{R}^d \times \mathbb{R}^d \rightarrow \mathbb{R} $ as
  \begin{equation}
    K_\tau(x_i,x_j) := \langle \partial_{\theta} f(x_j; \theta(\tau)) ,\, \partial_{\theta} f(x_i; \theta(\tau)) \rangle.
  \end{equation}
  This kernel is symmetric and positive semi-definite by construction.
  We remark that the NTK depends on the specific topology of the considered neural network and on the choice of the activation function.
  Considering wider neural networks, it can be proven \cite{du2018gradient, lee2019wide} that 
  \begin{equation}
  K_\tau(x_i, x_j) \underset{\mathbb{P}}{\rightarrow} K^{\infty}(x_i,x_j),
\end{equation}
where in the infinite-width limit $ K^{\infty}(x_i,x_j) \approx K_0(x_i, x_j)$, meaning that the average infinite-width limiting NTK is close to the NTK computed with weights and biases at the initialization.
Therefore, during the training processes of wide neural networks the trajectories of the cost functional ruled by $K_{\tau}$ are close to the approximated linearized ones (ruled by $K_\infty$).
\newpage
\section{KOL: algebraic derivation}
\label{app:kolDer}

We set the learning problem in the specific Reproducing Kernel Hilbert Spaces (RKHS) framework, so to employ the tools deriving from the Kernel Regression theory, thus enhancing the computational 
efficiency and capitalizing on the linearity of the resulting functional framework.
Let $\mathcal{U}$ be an RKHS of functions $u: \Omega \rightarrow \mathbb{R}$ where the kernel introduced is $Q: \Omega \times \Omega \rightarrow \mathbb{R}$, and $\mathcal{V}$ is an RKHS of functions $v: D \rightarrow \mathbb{R}$ endowed with the following kernel $K : D \times D \rightarrow \mathbb{R}$.
Then, $\psi$ and $\chi$ are defined as optimal recovery maps as
\begin{equation}
  \begin{split}
    \psi(U) &:= \argmin_{w \in \mathcal{U}} \| w \|_Q \; \; \mathrm{s.t.} \; \; \phi(w) = U,\\
    \chi(V) &:= \argmin_{w \in \mathcal{V}} \| w \|_K \; \; \mathrm{s.t.} \; \; \phi(w) = V,\\
  \end{split}
\end{equation}
with $\| w \|_Q = \sqrt{\langle \langle w, Q(\cdot,{t}) \rangle , \langle w, Q(\cdot,{t}) \rangle \rangle }$ and $\| w \|_K = \sqrt{\langle \langle w, K(\cdot,t) \rangle , \langle w, K(\cdot,t) \rangle \rangle }$ are the norms defined in the respective Hilbert spaces exploiting the kernel reproducibility property.
Assuming that $\phi$ and $\varphi$ are pointwise evaluations at specific collocation points as previously discussed, the optimal recovery maps have explicit closed forms deriving from kernel interpolation theory \cite{scholkopf2001generalized}: 
\begin{equation}
  \psi(U)({t}) = Q({t},{T})Q({T},{T})^{-1}U, \;\; \chi(V)(t) = K(t,T) K(T,T)^{-1}V,
  \label{eq:kerDef}
\end{equation}
where $Q({T},{T})$ and $K(T,T)$ are symmetric and definite positive matrices such that $Q({T},{T})_{ij} = Q({T}_i,{T}_j) $ and $K(T,T)_{ij} = K(T_i, T_j)$, whilst $Q({t},{T})_i = Q({t},{T}_i)$ and $K(t, T)_i = K(t, T_i)$ represent row-vectors.
Following the Kernel-based approach proposed in \cite{batlle2024kernel}, the operator learning scheme simplifies to determine an approximation of the mapping between two finite-dimensional Euclidean spaces $f^{\dagger}: \mathbb{R}^n \rightarrow \mathbb{R}^n$, defined as
\begin{equation}
  f^{\dagger} := \varphi \circ \mathcal{G} \circ \psi,
\end{equation}
cf. Figure \ref{fig:diagram}, where the reconstruction maps $\psi: \mathbb{R}^n \rightarrow \mathcal{U}$ and $\chi: \mathbb{R}^n \rightarrow \mathcal{V}$ need to be properly defined.

Finally, we aim at approximating the $f^{\dagger}$ function through the use of a vector-valued kernel.
Indeed, consider $\Gamma : \mathbb{R}^{n} \times \mathbb{R}^n \rightarrow \mathcal{L}(\mathbb{R}^n)$ the matrix valued kernel following the notation of \cite{alvarez2012kernels}.
We call the RKHS induced by the considered kernel $\mathcal{H}_{\Gamma}$, and the respective induced norm is $\| \cdot \|_{\mathcal{H}_{\Gamma}}$.
Hence, $f^{\dagger}$ can be approximated with the map $\bar{f}$ solving the following optimization problem:
\begin{equation}
  \bar{f} := \argmin_{f \in \mathcal{H}_{\Gamma}} \| f \|_{\Gamma} \; \; \mathrm{s.t.} \;\; f(\phi(u_i)) = \varphi(v_i), \; i = 1,2\hdots N.
\end{equation}

We introduce the following block vectors \textbf{U} and \textbf{V} as
\begin{equation}
  \textbf{U} = \begin{pmatrix}
    U_1\\
    U_2\\
    \vdots\\
    U_{N-1}\\
    U_N
  \end{pmatrix} \in \mathbb{R}^{n\,N}, \;   \textbf{V} = \begin{pmatrix}
    V_1\\
    V_2\\
    \vdots\\
    V_{N-1}\\
    V_N
  \end{pmatrix} \in \mathbb{R}^{n\,N},
\end{equation}
where $\{ U_i \}_i := \{ \phi(u_i) \}_i, \, \forall \, u_i \in \mathcal{U}$ and $\{ V_i \}_i := \{ \varphi(v_i) \}_i, \, \forall \, v_i \in \mathcal{V}$.
Then, we define with a slight abuse of notation the matrix $\Gamma: \mathbb{R}^{nN} \times \mathbb{R}^{nN} \rightarrow \mathbb{R}^{nN \times nN}$ as
\begin{equation}
  \Gamma(\textbf{U}, \textbf{U}) = \begin{bmatrix}
    \Gamma(U_1, U_1) & \Gamma(U_1, U_2) & \hdots & \Gamma(U_1, U_N)\\
    \Gamma(U_2, U_1) & \Gamma(U_2, U_2) & \hdots & \Gamma(U_2, U_N)\\
    \vdots & \vdots & & \vdots \\
    \Gamma(U_N, U_1) & \Gamma(U_N, U_2) & \hdots & \Gamma(U_N, U_N)\\
  \end{bmatrix},
  \label{eq:KUU}
\end{equation}
where each $\Gamma(U_i, U_j) \in \mathbb{R}^{n \times n}, \; \forall i,j=1,2 \hdots N $ is an independent block, and the following matrix
\begin{equation}
  \Gamma(U, \textbf{U}) = \begin{bmatrix}
    \Gamma(U, U_1) & \Gamma(U,U_2) & \hdots & \Gamma(U,U_N)
  \end{bmatrix} \in \mathbb{R}^{n \times nN}.
\end{equation}
In this work we assume to deal with uncorrelated input samples, therefore the analysis can be simplified by relying on diagonal kernels for $\Gamma$.
In particular, let $S : \mathbb{R}^n \times \mathbb{R}^n \rightarrow \mathbb{R}$ be a scalar kernel (see Section \ref{sec:kolEpi} where we compare different practical choices for $S$).
Thus, each block in \eqref{eq:KUU} is a diagonal block, \textit{i.e.}
\begin{equation}
  \Gamma(U_i, U_j) = S(U_i, U_j) I = \begin{bmatrix}
    S(U_i, U_j) & 0 & \hdots & 0\\
    0 & S(U_i, U_j) & \hdots & 0\\
    \vdots & \vdots & & \vdots \\
    0 & 0 & \hdots & S(U_i, U_j)\\
  \end{bmatrix}.
  \label{eq:GammaK}
\end{equation}

Hence, the problem of learning the operator from prescribed input-output couples can be recast as an optimal recovery problem:
\begin{equation}
  \bar{f}_j := \argmin_{g \in \mathcal{H}_S} \| g\|_{S} \; \mathrm{s.t.} \; g(\phi(u_i)) = \varphi(v_i)_j, \, \forall \, i = 1,2\hdots N, \, j= 1,2 \hdots n,
\end{equation}
where the RKHS endowed with the $S$ kernel is named $(\mathcal{H}_S, \| \cdot \|_S)$.
In this finite dimensional case, it is possible to employ the fundamental result from kernel theory known as representer theorem \cite{scholkopf2001generalized}.
Therefore, each component of the dicrete representation can be written explicitly as 
\begin{equation}
  \bar{f}_j(U) = S(U, \textbf{U}) S(\textbf{U}, \textbf{U})^{-1} \textbf{V}_{\cdot, j},
  \label{eq:zj}
\end{equation}
where $S(U, \textbf{U}): \mathbb{R}^n \times \mathbb{R}^{nN} \rightarrow \mathbb{R}^N $ is a row vector, $S(\textbf{U}, \textbf{U}): \mathbb{R}^{nN} \times \mathbb{R}^{nN} \rightarrow \mathbb{R}^{N \times N}$ and $\textbf{V}_{\cdot,j} = [[V_1]_j, [V_2]_j \hdots [V_N]_j]^T \in \mathbb{R}^N$.
Equation \eqref{eq:zj} in the alternative kernel-methods-like form 
\begin{equation}
  \bar{f}(U) =\sum_{j=1}^N \underbrace{S(U, U_j)}_{\in \mathbb{R}} \underbrace{\alpha_j}_{\in \mathbb{R}^n}.
  \label{eq:zj_2}
\end{equation}
Indeed, by the representer theorem it holds that
\begin{equation}
  V(U) = \Gamma(U, \textbf{U}) \Gamma(\textbf{U}, \textbf{U})^{-1} \textbf{V} = \Gamma(U, \textbf{U}) \boldsymbol{\alpha} = \sum_{j=1}^N \underbrace{\Gamma(U, U_j)}_{\in \mathbb{R}^{n \times n}} \alpha_j =\sum_{j=1}^N \underbrace{S(U, U_j)}_{\in \mathbb{R}} \alpha_j ,
  \label{eq:Vu}
\end{equation}
where 
\begin{equation}
  \boldsymbol{\alpha} = \begin{pmatrix}
    \alpha_1\\
    \alpha_2\\
    \vdots\\
    \underbrace{\alpha_N}_{\in \mathbb{R}^n}
  \end{pmatrix} =  {\tiny \begin{pmatrix}
  \begin{array}{c}
    \alpha_1^1\\
    \alpha_1^2\\
    \vdots\\
    \alpha_1^n\\
    \hdashline\\
    \alpha_2^1\\
    \alpha_2^2\\
    \vdots\\
    \alpha_2^n\\
    \hdashline\\
    \vdots\\
    \hdashline\\
    \alpha_N^1\\
    \alpha_N^2\\
    \vdots\\
    \alpha_N^n\\
    \end{array}
  \end{pmatrix}}\in \mathbb{R}^{nN} \; \mathrm{s.t} \; \Gamma(\textbf{U}, \textbf{U}) \boldsymbol{\alpha} = \textbf{V}.
  \label{eq:initialKernel}
\end{equation}
The matrix in system in \eqref{eq:initialKernel} can be written explicitly as
\begin{equation}
\begin{split}
  &\Gamma(\textbf{U}, \textbf{U}) = \\
 &\begin{bmatrix} \tiny
    \begin{array}{c@{\hspace{-10pt}}c@{\hspace{3pt}}c@{\hspace{2pt}}c:c@{\hspace{-10pt}}c@{\hspace{3pt}}c@{\hspace{2pt}}c:c:c@{\hspace{-10pt}}c@{\hspace{3pt}}c@{\hspace{2pt}}c}    
       S(U_1,U_1) &  0 &  \hdots &  0 & S(U_1,U_2) &  0 &  \hdots &  0 & \hdots & S(U_1,U_N) &  0 &  \hdots &  0\\
      0 & S(U_1,U_1) &  \hdots &  0 & 0 & S(U_1,U_2) &  \hdots &  0 & \hdots & 0 & S(U_1,U_N) &  \hdots &  0\\
      \vdots & & & \vdots & \vdots & & & \vdots & \hdots & \vdots & & & \vdots\\
      0 & 0 &  \hdots &  S(U_1,U_1) & 0 & 0 &  \hdots &  S(U_1,U_2) & \hdots & 0 & 0 &  \hdots &  S(U_1,U_N) \\
      \hdashline
      S(U_2,U_1) &  0 &  \hdots &  0 & S(U_2,U_2) &  0 &  \hdots &  0 & \hdots & S(U_2,U_N) &  0 &  \hdots &  0\\
      0 & S(U_2,U_1) &  \hdots &  0 & 0 & S(U_2,U_2) &  \hdots &  0 & \hdots & 0 & S(U_2,U_N) &  \hdots &  0\\
      \vdots & & & \vdots & \vdots & & & \vdots & \hdots & \vdots & & & \vdots\\
       0 & 0 &  \hdots &  S(U_2,U_1) & 0 & 0 &  \hdots &  S(U_2,U_2) & \hdots & 0 & 0 &  \hdots &  S(U_2,U_N) \\
       \hdashline
       \vdots & & & \vdots & \vdots & & & \vdots & \hdots & \vdots & & & \vdots\\  
       \hdashline
       S(U_N,U_1) &  0 &  \hdots &  0 & S(U_N,U_2) &  0 &  \hdots &  0 & \hdots & S(U_N,U_N) &  0 &  \hdots &  0\\
      0 & S(U_N,U_1) &  \hdots &  0 & 0 & S(U_N,U_2) &  \hdots &  0 & \hdots & 0 & S(U_N,U_N) &  \hdots &  0\\
      \vdots & & & \vdots & \vdots & & & \vdots & \hdots & \vdots & & & \vdots\\
       0 & 0 &  \hdots &  S(U_N,U_1) & 0 & 0 &  \hdots &  S(U_N,U_2) & \hdots & 0 & 0 &  \hdots &  S(U_N,U_N) \\
    \end{array}
  \end{bmatrix}
    \end{split}
\end{equation}

Hence, we can reorder the system as 
\begin{equation}
  \begin{cases}
    S(U_1, U_1) \alpha_1^1 + S(U_1,U_2) \alpha_2^1 + S(U_1,U_3) \alpha_3^1 + \hdots + S(U_1,U_N) \alpha_N^1 = V_1^1\\
    S(U_2, U_1) \alpha_1^1 + S(U_2,U_2) \alpha_2^1 + S(U_2,U_3) \alpha_3^1 + \hdots + S(U_2,U_N) \alpha_N^1 = V_2^1\\
    S(U_1, U_1) \alpha_1^1 + S(U_1,U_2) \alpha_2^1 + S(U_1,U_3) \alpha_3^1 + \hdots + S(U_1,U_N) \alpha_N^1 = V_1^1\\
    S(U_3, U_1) \alpha_1^1 + S(U_3,U_2) \alpha_2^1 + S(U_3,U_3) \alpha_3^1 + \hdots + S(U_3,U_N) \alpha_N^1 = V_3^1\\ 
    \vdots\\
    S(U_N, U_1) \alpha_1^1 + S(U_N,U_2) \alpha_2^1 + S(U_N,U_3) \alpha_3^1 + \hdots + S(U_N,U_N) \alpha_N^1 = V_N^1\\     
  \end{cases},
\end{equation}
and we can define the matrix
\begin{equation}
  \textbf{S}(\textbf{U}, \textbf{U}) = \begin{bmatrix}
    S(U_1, U_1) & S(U_1,U_2) & S(U_1,U_3) & \hdots & S(U_1, U_N)\\
    S(U_2, U_1) & S(U_2,U_2) & S(U_2,U_3) & \hdots & S(U_2, U_N)\\
    \vdots & & & & \vdots\\
    S(U_N, U_1) & S(U_N,U_2) & S(U_N,U_3) & \hdots & S(U_N, U_N)\\    
  \end{bmatrix} \in \mathbb{R}^{N \times N}. 
\end{equation}
We get that 
\begin{equation}
  \alpha_i = \begin{pmatrix}
    [\textbf{S}(\textbf{U}, \textbf{U})^{-1} \textbf{V}_{\cdot, 1}]_i\\
    [\textbf{S}(\textbf{U}, \textbf{U})^{-1} \textbf{V}_{\cdot, 2}]_i\\
    \vdots\\
    [\textbf{S}(\textbf{U}, \textbf{U})^{-1} \textbf{V}_{\cdot, n}]_i\\  
  \end{pmatrix}, \; \forall i = 1,2\hdots N.
\end{equation}
Finally, we rewrite equation \eqref{eq:Vu} as
\begin{equation}
  \bar{f}_j(U) = V(U)_j = \left [ \sum_{i=1}^{N} S(U,U_i) \alpha_i \right ]_j = \left [ \sum_{i=1}^{N} S(U,U_i) \begin{pmatrix}
  \begin{array}{c}
    [\textbf{S}(\textbf{U}, \textbf{U})^{-1} \textbf{V}_{\cdot, 1}]_i\\
    [\textbf{S}(\textbf{U}, \textbf{U})^{-1} \textbf{V}_{\cdot, 2}]_i\\
    \vdots\\
    \hdashline\\
    [\textbf{S}(\textbf{U}, \textbf{U})^{-1} \textbf{V}_{\cdot, j}]_i\\
    \hdashline
    \vdots\\
    [\textbf{S}(\textbf{U}, \textbf{U})^{-1} \textbf{V}_{\cdot, n}]_i\\
    \end{array}
  \end{pmatrix} \right ]_j,
\end{equation}
that, by rearranging the terms is exactly equation \eqref{eq:zj}.

Combining equations \eqref{eq:kerDef} and \eqref{eq:zj} we obtain the approximation operator as 
\begin{equation}
  \bar{\mathcal{G}} := \chi \circ \bar{f} \circ \phi.
\end{equation}
Finally, recalling the definition of the operator in equation \eqref{eq:zj_2}, $\bar{\mathcal{G}}: \mathcal{U} \rightarrow \mathcal{V}$ is 
\begin{equation}
  \bar{\mathcal{G}}(u) = \chi \left ( \sum_{j=1}^N S(\phi(u), U_j) \alpha_j \right ).
\end{equation}
\end{document}